\documentclass[12pt]{article}
\usepackage{amssymb,amsmath,bm}
\usepackage{mathrsfs}
\usepackage{hyperref}
\hypersetup{
    colorlinks,
    citecolor=black,
    filecolor=black,
    linkcolor=black,
    urlcolor=black
}

\usepackage{color}

\newcommand{\elaw}{\stackrel{\rm Law}{=}}
\newcommand{\claw}{\stackrel{\rm Law}{\longrightarrow}}
\newcommand{\ctv}{\stackrel{TV}{\longrightarrow}}

\usepackage{authblk}

\usepackage{constants}
\usepackage{enumerate}
\begin{document}
\newcommand{\beq}{\begin{eqnarray}}
\newcommand{\eeq}{\end{eqnarray}}
\newcommand{\beas}{\begin{eqnarray*}}
\newcommand{\enas}{\end{eqnarray*}}
\newcommand{\bea}{\begin{eqnarray}}
\newcommand{\ena}{\end{eqnarray}}
\newcommand{\bms}{\begin{multline*}}
\newcommand{\ems}{\end{multline*}}
\newcommand{\qmq}[1]{\quad \mbox{#1} \quad}
\newcommand{\qm}[1]{\quad \mbox{#1}}
\newcommand{\nn}{\nonumber}
\newcommand{\bbox}{\hfill $\Box$}
\newcommand{\ignore}[1]{}
\newcommand{\bs}[1]{\boldsymbol{#1}}
\newcommand{\tr}{\mbox{tr}}
\newtheorem{theorem}{Theorem}[section]
\newtheorem{corollary}{Corollary}[section]
\newtheorem{conjecture}{Conjecture}[section]
\newtheorem{proposition}{Proposition}[section]
\newtheorem{remark}{Remark}[section]
\newtheorem{lemma}{Lemma}[section]
\newtheorem{definition}{Definition}[section]
\newtheorem{condition}{Condition}[section]
\newtheorem{example}{Example}[section]
\newcommand{\pf}{\noindent {\bf Proof:} }
\def\blfootnote{\xdef\@thefnmark{}\@footnotetext}

\title{{ \bf\Large Gaussian Phase Transitions and Conic Intrinsic Volumes:  Steining the Steiner Formula}}
\author[1]{Larry Goldstein}
\author[2]{Ivan Nourdin}
\author[2]{Giovanni Peccati}
\affil[1]{University of Southern California}
\affil[2]{University of Luxembourg}
 \maketitle
\date{}

\begin{abstract}
Intrinsic volumes of convex sets are natural geometric quantities that also play important roles in applications, such as linear inverse problems with convex constraints,  and constrained statistical inference. It is a well-known fact that, given a closed convex cone $C\subset \mathbb{R}^d$, its conic intrinsic volumes determine a probability measure on the finite set $\{0,1,...d\}$, customarily denoted by $\mathcal{L}(V_C)$. The aim of the present paper is to {provide a Berry-Esseen bound for the normal approximation of ${\cal L}(V_C)$, implying} a general quantitative central limit theorem (CLT) for sequences of (correctly normalised) discrete probability measures of the type $\mathcal{L}(V_{C_n})$, $n\geq 1$. This bound shows that, in the high-dimensional limit, most conic intrinsic volumes encountered in applications can be approximated by a suitable Gaussian distribution. Our approach is based on a variety of techniques, namely: (1) Steiner formulae for closed convex cones, (2) Stein's method and second order Poincar\'e inequality, (3) concentration estimates, and (4) Fourier analysis. { Our results explicitly connect the sharp phase transitions, observed in many regularised linear inverse problems with convex constraints, with the asymptotic Gaussian fluctuations of the intrinsic volumes of the associated descent cones. In particular, our findings complete and further illuminate the recent breakthrough discoveries by Amelunxen, Lotz, McCoy and Tropp (2014) and McCoy and Tropp (2014) about the concentration of conic intrinsic volumes and its connection with threshold phenomena}. As an additional outgrowth of our work we develop total variation bounds for normal approximations of the lengths of projections of Gaussian vectors on {closed} convex sets. \end{abstract}

\section{Introduction}

\subsection{Overview}

Every closed convex cone $C \subset \mathbb{R}^d$ can be associated with a random variable $V_C$, with support on $\{0,\ldots,d\}$ {whose distribution ${\cal L}(V_C)$ coincides} with the so-called {\it conic intrinsic volumes} of $C$. The distribution ${\cal L}(V_C)$ is a natural object that summarizes key information about the geometry of $C$, and is important in applications, ranging from compressed sensing to constrained statistical inference. In particular, for a closed convex cone $C$ the mean $\delta_C=EV_C$ (which is customarily called the {\it statistical dimension} of $C$) measures in some sense the `effective' dimension of $C$, and generalises the classical notion of dimension for linear subspaces. As proved in the groundbreaking papers Amelunxen, Lotz, McCoy and Tropp \cite{AmLo13} and by McCoy and Tropp \cite{McTr13} (see also Section \ref{ss:app1} below for a more detailed discussion of this point), in the case of the so-called {\it descent cones} arising in convex optimisation, the concentration of the distribution of $V_C$  around $\delta_C$ explains with striking precision threshold phenomena exhibited by the probability of success in linear inverse problems with convex constraints.

\medskip

Our principal aim in this paper is {to produce a Berry-Esseen bound for ${\cal L}(V_C)$ leading to}  minimal conditions on a sequence  of closed convex cones $\{C_n\}_{n\geq 1}$, ensuring that the sequence
\beas
\frac{V_{C_n}-EV_{C_n}}{\sqrt{{\rm Var}(V_{C_n})}}, \quad n\geq 1,
\enas
converges in distribution towards a standard Gaussian ${\cal N}(0,1)$ random variable. The bounds in our main findings depend only on the mean and the variance of the random variables $V_{C_n}$, and are summarized in Part 2 of Theorem \ref{t:mainabstract} below.

\medskip

As explained in the sections to follow, the strategy for achieving our goals consists in using the elegant {\it Master Steiner formula} from McCoy and Tropp \cite{McTr13}, in order to connect random variables of the type $V_C$ to objects with the form $\|\Pi_C({\bf g})\|^2$, where ${\bf g}$ is a standard Gaussian vector, $\Pi_C$ is the metric projection onto $C$, and $\|\cdot\|$ stands for the Euclidean norm. Shifting from $V_C$ to $\|\Pi_C({\bf g})\|^2$ allows one to unleash the full power of some recently developed techniques for normal approximations, based on the interaction between Stein's method (see \cite{ChGoSh10}) and variational analysis on a Gaussian space (see \cite{NoPe12}). In particular, our main tool will be the so-called {\it second order Poincar\'e inequality} developed in \cite{Ch09, NoPeRe09}. In Section 4, we will also use techniques from Fourier analysis in order to compute explicit Berry-Esseen bounds.

\medskip

As discussed below, our findings represent a significant extension of the results of \cite{AmLo13, McTr13}, where the concentration of ${\cal L}(V_C)$ around $\delta_C$ was first studied by means of tools from Gaussian analysis, as well as by exploiting the connection between intrinsic volumes and metric projections. { Explicit applications to regularised linear inverse problems are described in detail in Section \ref{ss:app1} below}.

\medskip

We will now quickly present some basic facts of conic geometry that are relevant for our analysis. Our main theoretical contributions are discussed in Section \ref{ss:main}, whereas connections with applications are described in Section \ref{ss:app1} and Section \ref{ss:app2}.

\subsection{Elements of conic geometry}\label{ss:intro conic}

The reader is referred to the classical references \cite{Roc, rw}, as well as to \cite{AmLo13, McTr13}, for any unexplained notion or result related to convex analysis.

\medskip

\noindent\underline{\it Distance from a convex set and metric projections.} Fix an integer $d\geq 1$. Throughout the paper, we shall denote by $\langle {\bf x}, {\bf y} \rangle$ and $\|{\bf x}\|^2=\langle {\bf x}, {\bf x} \rangle$, respectively, the standard inner product and squared Euclidean norm in $\mathbb{R}^d$. Given a closed convex set $C\subset \mathbb{R}^d$, we define the {\it distance} between a point ${\bf x}$ and $C$ as
\bea \label{def:dxC}
d({\bf x},C):=\inf_{{\bf y} \in C}\|{\bf x}-{\bf y}\|.
\ena
By the strict convexity of the mapping ${\bf x}\mapsto \|{\bf x}\|^2$, the infimum is attained at a unique vector, called the {\it metric projection} of ${\bf x}$ onto $C$, which we denote by $\Pi_C({\bf x})$.
\medskip

\noindent\underline{\it Convex cones and polar cones.} A set $C \subset \mathbb{R}^d$ is a {\it convex cone} if $a {\bf x} + b {\bf y} \in C$ whenever ${\bf x}$ and ${\bf y}$ are in $C$ and $a$ and $b$ are positive reals. The {\it polar cone} $C^0$ of a cone $C$ is given by
\bea \label{def:polar.bear}
C^0 = \left\{{\bf y} \in \mathbb{R}^d: \langle {\bf y},{\bf x} \rangle \le 0, \forall {\bf x} \in C \right\}.
\ena
It is easy to verify that the polar cone of a closed convex cone is again a closed convex cone. By virtue e.g. of \cite[formula (7.2)]{McTr13}, 
any vector ${\bf x} \in \mathbb{R}^d$ may be written as:
\bea \label{eq:polar.orth.decomp}
{\bf x} = \Pi_C({\bf x}) + \Pi_{C^0}({\bf x}) \qmq{with $\Pi_C({\bf x}) \perp \Pi_{C^0}({\bf x}),$}
\ena
where the orthogonality relation is in the sense of the inner product $\langle \cdot, \cdot \rangle$ on $\mathbb{R}^d$. A quick computation shows also that, for every closed convex cone $C$ and every ${\bf x}\in \mathbb{R}^d$,
\bea
 \|\Pi_C({\bf x})\| = \sup_{{\bf y}\in C : \|{\bf y}\|\leq 1} \langle {\bf x}, {\bf y}\rangle.\label{e:bobo}
\ena

\medskip

\noindent\underline{\it Steiner formulae and intrinsic volumes.} Letting $B^d$ and $S^{d-1}$ denote, respectively, the unit ball and unit sphere in $\mathbb{R}^d$, the classical {\it Steiner formula} for the Euclidean expansion of a compact convex set $K$ states that
\beas
\mbox{Vol}(K+\lambda B_d) = \sum_{j=0}^d \lambda^{d-j}\mbox{Vol}(B^{d-j}) {\mathcal V}_j \qmq{for all $\lambda \ge 0$,}
\enas
where addition on the left-hand side indicates the Minkowski sum of sets, and the numbers ${\mathcal V}_j, j=0,\ldots,d$ on the right, called {\it Euclidean intrinsic volumes}, depend only on $K$. The Euclidean intrinsic volumes numerically encode key geometric properties of $K$, for instance, ${\mathcal V}_d$ is the volume, $2{\mathcal V}_{d-1}$ the surface area, and ${\mathcal V}_0$ the Euler characteristic of $K$. See e.g. \cite[p. 142]{at}, \cite[Chapter 7]{kr} and \cite[p. 600]{sw} for standard proofs.

An `angular' Steiner formula was developed in \cite{Al48, He43, Sa50}, and expresses the size of an angular expansion of a closed convex cone $C$ as follows:
\bea \label{angular.Steiner}
P\left\{ d^2({\bs \theta},C) \le \lambda \right\}= \sum_{j=0}^d \beta_{j,d}(\lambda) v_j,
\ena
where ${\bs \theta}$ is a random variable uniformly distributed on $S^{d-1}$, the coefficients
$$
\beta_{j,d} (\lambda)= P[B(d-j,d)\leq \lambda]
$$
(where each $B(d-j,d)$ has the Beta distribution with parameters $(d-j)/2$ and $d/2$) do not depend on $C$, and the {\it conic intrinsic volumes} $v_0,\ldots,v_d$ are determined by $C$ only, and can be shown to be nonnegative and sum to one. As a consequence, we may associate to the conic intrinsic volumes of $C$ an integer-valued random variable $V$, whose probability distribution ${\cal L}(V)$ is given by
\begin{equation}\label{VC}
P(V=j)=v_j, \qmq{for $j=0,\ldots,d$.}
\end{equation}
When the dependence of any quantities on the cone needs to be emphasized, we will write $V_C$ for $V$ and $v_j(C)$ for $v_j$, $j=0,\ldots,d$. As shown in \cite{McTr13}, relation \eqref{angular.Steiner} can be seen as a consequence of a general result, known as {\it Master Steiner formula} and stated formally in Theorem \ref{t:msf} below. Such a result implies that, writing ${\bf g}\sim {\cal N}(0,I_d)$ for a standard $d$-dimensional Gaussian vector, the squared norms $\|\Pi_C({\bf g})\|^2$ and $\|\Pi_{C^0}({\bf g})\|^2$ behave like two independent chi-squared random variables with a random number $V_C$ and $d-V_C$, respectively, of degrees of freedom: in symbols,
\bea\label{e:rita}
( \|\Pi_C({\bf g})\|^2, \|\Pi_{C^0}({\bf g})\|^2 ) \sim (\chi^2_{V_C},\chi^2_{d-V_C}).
\ena
In particular, equation \eqref{e:rita} is consistent with the well-known relation $v_j(C) = v_{d-j}(C^0)$ ($j=0,...,d$), that is: the distribution of the random variable $V_{C^0}$, associated with the polar cone $C^0$ via its intrinsic volumes, satisfies the relation
\bea \label{eq:VC.polar.bear}
V_{C^0} \elaw d-V_C,
\ena
where, here and in what follows, $\elaw$ indicates equality in distribution. To conclude, we notice that partial versions of \eqref{e:rita} (only involving $\|\Pi_C({\bf g})\|^2$) were already known in the literature prior to \cite{McTr13}, in particular in the context of constrained statistical inference --- see e.g. \cite{d, sh0, sh}, as well as \cite[Chapter 3]{ss}.

\medskip

\noindent\underline{\it Statistical dimensions}. As for Euclidean intrinsic volumes, the distribution of $V_C$ encodes key geometric properties of $C$. For instance, the mean $\delta_C:= E[V_C] =E\|\Pi_C({\bf g})\|^2$, generalizes the notion of dimension. In particular, if $L_k$ is a linear  subspace of $\mathbb{R}^d$ of dimension $k$, and hence a closed convex cone, then $v_j(L_k)$ is one when $j=k$ and zero otherwise, and therefore $\delta(L_k)=k$. The parameter $\delta_C$ is often called the {\it statistical dimension} of $C$. We observe that, in view of \eqref{e:bobo}, the statistical dimension $\delta_C$ is tightly related to the so-called {\it Gaussian width} of a convex cone
$$
w_C := E\left(\sup_{{\bf y}\in C : \|{\bf y}\|\leq 1} \langle {\bf g}, {\bf y}\rangle\right),
$$
where ${\bf g}\sim {\cal N}(0,I_d)$. The notion of Gaussian width plays an important role in many key results of compressed sensing (see e.g. \cite{RuVe08}). Standard arguments yield that $w^2_C\leq \delta_C\leq w^2_C+1$ (see \cite[Proposition 10.2]{AmLo13}). One situation where the statistical dimension is particularly simple to calculate is when $C$ is self dual, that is, when $C=-C^0$. In this case, $\delta_C=d/2$ by \eqref{eq:VC.polar.bear}. The nonnegative orthant, the
second-order cone, and the cone of positive-semidefinite matrices are all self dual; see \cite{McTr13} for definitions and further explanations.

\medskip

\noindent\underline{\it Polyhedral cones.} We recall that a {\it polyhedral cone} $C$ is one that can be expressed as the intersection of a finite number of halfspaces, that is, one for which there exists an integer $N$ and vectors ${\bf u}_1,\ldots,{\bf u}_N$ in $\mathbb{R}^d$ such that
\beas
C= \bigcap_{i=1}^N \{{\bf x} \in \mathbb{R}^d: \langle {\bf u}_i, {\bf x} \rangle \ge 0\}.
\enas
For polyhedral cones the probabilities $v_j, j=0,\ldots,d$ can be connected to the behavior of the projection $\Pi_C({\bf g})$ of a standard Gaussian variable ${\bf g} \sim {\cal N}(0,I_d)$ onto $C$. Indeed, in this case we have the representation
\bea \label{proj.relint.j}
v_j=P\left( \Pi_C({\bf g}) \,\mbox{lies in the relative interior of a $j$-dimensional face of $C$}\right)
\ena
(see e.g. \cite{AmLo13, McTr13}).

\subsection{Main theoretical contributions}\label{ss:main}

The main result of the present paper is the following general central limit theorem (CLT), involving the intrinsic volume distributions of a sequence of closed convex cones with increasing statistical dimensions.

\begin{theorem}\label{t:mainabstract} Let $\{d_n:n\ge 1\}$ be a sequence of non-negative integers and let $\{C_n \subset \mathbb{R}^{d_n}: n \ge 1\}$ be a collection of non-empty closed convex cones such that $\delta_{C_{n}} \rightarrow \infty$, and write $\tau^2_{C_n} = {\rm Var}(V_{C_n})$, $n\geq 1$. For every $n$, let ${\bf g}_n\sim {\cal N}(0,I_{d_n})$ and write $\sigma^2_{C_n} = {\rm Var}(\|\Pi_{C_n}({\bf  g}_n)\|^2)$, $n\geq 1$. Then, the following holds.
\begin{enumerate}
\item[\rm 1.] One has that $2\delta_{C_n} \leq \sigma_{C_n}^2\leq 4\delta_{C_n}$ for every $n$ and, as $n\to\infty$, the sequence
$$
\frac{\|\Pi_{C_n}({\bf  g}_n)\|^2 - \delta_{C_n}}{\sigma_{C_n}}, \quad n\geq 1,
$$
converges in distribution to a standard Gaussian random variable $N\sim \mathcal{N}(0,1)$.

\item[\rm 2.] If, in addition, $\liminf_{n\to\infty} \tau^2_{C_n}/\delta_{C_n} >0$, then, as $n\to\infty$, the sequence
$$
\frac{V_{C_n} - \delta_{C_n}}{\tau_{C_n}}, \quad n\geq 1,
$$
also converges in distribution to $N\sim \mathcal{N}(0,1)$, and {moreover one has} the Berry-Esseen estimate
\begin{equation}\label{e:joe}
\sup_{u\in \mathbb{R}}\left| P\left[\frac{V_{C_n} - \delta_{C_n}}{\tau_{C_n}}\leq u\right] - P[N\leq u]\right| = O\left(\frac{1}{\sqrt{\log\delta_{C_n}}}\right).
\end{equation}
\end{enumerate}
\end{theorem}

Part 1 follows from Corollary \ref{cor:deltabound}. Part 2 is a consequence of Theorem \ref{thm:Linfty.VC} 
below that provides a Berry-Esseen bound, with small explicit constants, for { the normal approximation of $V_C$ and} for any closed convex cone $C$, in terms of $\delta_C,\sigma_C^2$ and $\tau_C^2$. { In particular, if $C$ is a closed convex cone such that $\tau_C>0$, then we will prove in Theorem \ref{thm:Linfty.VC}
and Remark \ref{rmk:Linfty.VC} that, writing $\alpha := \tau_C^2/\delta_C$, for $\delta_C \ge 8$, 
\begin{equation}\label{e:kintro}
\sup_{u\in \mathbb{R}}\left| P\left[\frac{V_{C} - \delta_{C}}{\tau_{C}}\leq u\right] - P[N\leq u]\right|\leq h(\delta_C) + \frac{48}{\sqrt{\alpha \log^+(\alpha\sqrt{2}\delta_C)}},
\end{equation}
where
\begin{equation}\label{e:xox}
h(\delta) =\frac{1}{72} \left(\frac{\log \delta }{\delta^{3/16}}\right)^{5/2}.
\end{equation}

\begin{remark}{\rm Observe that, if one considers the sequence $\{C_d\}_{d \geq 1}$ consisting of the non-negative orthants of $\mathbb{R}^d$ , then $V_{C_d}$ follows a binomial distribution with parameters $(1/2, d)$ (in particular, $\delta_{C_d} = d/2$). It follows that, in this case, the supremum on the left-hand side of \eqref{e:joe} converges to zero at a speed of the order $O(d^{-1/2})$, from which we conclude that the rate supplied by \eqref{e:joe} is, in general, not optimal. }
\end{remark}

As anticipated, our strategy for proving Theorem \ref{t:mainabstract} (exception made for the Berry-Esseen bound \eqref{e:joe}) is to connect the distributions of $\|\Pi_{C_n}({\bf  g}_n)\|^2$ and $V_{C_n}$ via the Master Steiner formula \eqref{e:rita}, and then to study the normal approximation of the squared norm of $\Pi_{C_n}({\bf  g}_n)$ by means of Stein's method, as well as of general variational techniques on a Gaussian space (see \cite{ChGoSh10, NoPe12}). As illustrated in the Appendix contained in Section 5 below, Stein's method proceeds by manipulating a characterizing equation for a target distribution (in this case the normal), typically through couplings or integration by parts. Hence, we justify the title of this work by the heavy use that our application of Stein's method makes of relation \eqref{e:rita}, generalizing the angular Steiner formula \eqref{angular.Steiner}. As mentioned above, our main tool will be a form of the second order Poincar\'e inequalities studied in \cite{Ch09, NoPeRe09}.

{
\begin{remark}{\rm A crucial point one needs to address when applying Part 2 of Theorem \ref{t:mainabstract} is that, in order to check the assumption $\liminf_{n\to\infty} \tau^2_{C_n}/\delta_{C_n} >0$, one has to produce an effective lower bound on the sequence of conic variances $\tau^2_{C_n}$, $n\geq 1$. This issue is dealt with in Section \ref{s:ulbounds}, where we will prove new upper and lower bounds for conic variances, by using an improved version of the Poincar\'e inequality (see Theorem \ref{thm:Poincare}), as well as a representation of the covariance of smooth functionals of Gaussian fields in terms of the Ornstein-Uhlenbeck semigroup, as stated in formula \eqref{Poinc.Covariance} below. In particular, our main findings of Section \ref{s:ulbounds} (see Theorem \ref{t:lowvc}) will indicate that, in many crucial examples, the sequence $n\mapsto \tau^2_{C_n}$ {eventually satisfies} a relation of the type $$ c \| E[\Pi_{C_n}({\bf g})] \|^2 \leq \tau^2_{C_n} \leq 2 \| E[\Pi_{C_n}({\bf g})] \|^2,$$ where $c\in (0,2)$ does not depend on $n$. In view of Jensen inequality, this conclusion strictly improves the estimate $\tau_{C_n}^2\leq 2 \delta_{C_n}$ that one can derive e.g. from \cite[Theorem 4.5]{McTr13}.
}
\end{remark}}

We obtain normal approximation results for random variables that are more general than $\|\Pi_{C}({\bf  g})\|^2$. To this end, fix a closed convex cone $C\subset \mathbb{R}^d$ and ${\bs \mu}\in \mathbb{R}^d$, and introduce the shorthand notation:
\bea \label{def:F.sigma2.cone.gen}
F=\|\bs{\mu}-\Pi_C({\bf g}+\bs{\mu})\|^2 - m, \qmq{with} m=E[\|\bs{\mu}-\Pi_C({\bf g}+\bs{\mu})\|^2] \qmq{and} \sigma^2={\rm Var}(F).
\ena
Then, we prove in Theorem \ref{thm:main.cone} that
\bea\label{e:hairjoke}
 d_{{TV}}(F,N) \le \frac{16}{\sigma^2}\left\{ { \sqrt{m}(1+2\|{\bs \mu}\|) + 3\|{\bs \mu}\|^2+\|{\bs \mu}\| } \right\},
\ena
where $N \sim {\cal N}(0,\sigma^2)$ and $d_{TV}$ stand for the total variation distance, defined in \eqref{def:dtv.sets},  between the distribution of two random variables. In the fundamental case ${\bs \mu} ={\bf 0}$, Proposition \ref{cor:deltabound} shows that the previous estimate implies the simple relation
\bea\label{e:hairjoke2}
d_{{TV}}\left(\|\Pi_C({\bf g})\|^2-\delta_C,\,N\right) \le \frac{8}{\sqrt{\delta_C}},
\ena
where $N \sim {\cal N}(0,\sigma_C^2)$. Relation \eqref{e:hairjoke2} reinforces our intuition that the statistical dimension $\delta_C$ encodes a crucial amount of information about the distributions of $\|\Pi_C({\bf g})\|^2$ and, therefore, {about} $V_C$, via \eqref{e:rita}.

\medskip

It does not seem possible to directly combine the powerful inequality \eqref{e:hairjoke2} with \eqref{e:rita} in order to deduce an explicit Berry-Esseen bound such as \eqref{e:joe}. This estimate
is obtained in Section 5, by means of Fourier theoretical arguments of a completely different nature.

\begin{remark}{\rm We stress that the crucial idea that one can study a random variable of the type $V_C$, by applying techniques of Gaussian analysis to the associated squared norm $\|\Pi_C({\bf g})\|^2$, originates from the path-breaking references \cite{AmLo13, McTr13}, where this connection is exploited in order to obtain explicit concentration estimates via the entropy method,  see \cite{Bo98} and \cite{Le05}.
}
\end{remark}

\medskip

As stated in the Introduction, we will now show that our results can be used to exactly characterise phase transitions in regularised inverse problems with convex constraints.

\subsection{Applications to exact recovery of structured unknowns}\label{ss:app1}

\subsubsection{General framework}

{ In what follows,} we give a summary of how the conic intrinsic volume distribution plays a role in convex optimization for the recovery of structured unknowns and refer the reader e.g. to the excellent discussions { in  \cite{AmLo13, c, ChRe12,  McTr13}} for more detailed information.

\medskip

In certain high dimension recovery problems some small number of observations may be taken on an unknown high dimensional vector or matrix ${\bf x}_0$, thus determining that the unknown lies in the feasible set ${\cal F}$ of all elements consistent with what has been observed. As ${\cal F}$ may be large, the recovery of ${\bf x}_0$ is not possible without additional assumptions, such as that the unknown possesses some additional structure such as being sparse, or of low rank. As searching ${\mathcal F}$ for elements possessing the given structure can be computationally expensive, one instead may consider a convex optimization problem of finding ${\bf x} \in \mathcal{F}$ that minimizes $f({\bf x})$ for some
proper convex function\footnote{a convex function having at least one finite value and never taking the value $-\infty$} that promotes the structure desired.

\medskip

The analysis of such an optimization procedure leads one naturally to the study of the {\em descent cone} ${\cal D}(f,{\bf x})$ of $f$ at the point ${\bf x}$, given by
\beas
{\cal D}(f, {\bf x}) = \{ {\bf y}: \exists \tau>0\mbox{ such that }f({\bf x}+\tau {\bf y}) \le f({\bf x}) \}.
\enas
That is, ${\cal D}(f,{\bf x})$ is the conic hull of all directions that do not increase $f$ near ${\bf x}$.
The proof of Part 1 of Theorem \ref{flag} below -- included here for completeness -- reflects the general result,
that in the case where ${\cal F}$ is a subspace, the convex optimization just described successfully recovers the unknown ${\bf x}_0$ if and only if
\bea \label{eq:twin.bad.hair.day}
{\cal F} \cap ({\bf x}_0+{\cal D}(f,{\bf x}_0))=\{{\bf x}_0\}
\ena
(see Section 4 of \cite{RuVe08} and Proposition 2.1 \cite{ChRe12}, and Fact 2.8 of \cite{AmLo13}).

\medskip

The work \cite{ChRe12} provides a systematic way { according to which} an appropriate convex function $f$ may be chosen to promote a given structure. When an unknown vector, or matrix, is expressed as a linear combination
\bea \label{eq:x.sum.atoms}
{\bf x}_0=c_1{\bf a_1} + \cdots c_k {\bf a}_k
\ena
for $c_i \ge 0$, $a_i \in {\mathcal A}$ a set of building blocks or atoms of vectors or matrices, and $k$ small, then one minimizes
\bea \label{eq:f.from.A}
f({\bf x})=\inf\{t>0: {\bf x} \in t {\rm conv}({\mathcal A}) \},
\ena
over the feasible set, where ${\rm conv}({\mathcal A})$ is the convex hull of ${\mathcal A}$.

\subsubsection{Recovery of sparse vectors via $\ell_1$ norm minimization}

{ We now consider}
the underdetermined linear inverse problem of recovering a sparse vector ${\bf x}_0 \in \mathbb{R}^d$
from the observation of ${\bf z}=A{\bf x}_0$, where for $m<d$ the known matrix $A \in \mathbb{R}^{m \times d}$ has independent entries each with the standard normal ${\cal N}(0,1)$ distribution. {We say the vector ${\bf x}_0$ is $s$-sparse if it has exactly $s$ nonzero components; the value of $s$ is typically much smaller than $d$.} As a sparse vector is a linear combination of a small number of standard basis vectors, the prescription \eqref{eq:f.from.A} leads us to find a feasible vector that minimizes the $\ell_1$ norm, denoted by $\|\cdot \|_1$. { It is a well-known fact that such a linear inverse problem displays a {\it sharp phase transition} (sometimes called a {\it threshold phenomenon}): heuristically, this means that, for every value of $d$, there exists a very narrow band $[m_1, m_2]$ (that depends on $d$ and on the sparsity level of ${\bf x}_0$) such that the probability of {recovering ${\bf x}_0$ exactly is negligible} for $m<m_1$, and overwhelming for $m>m_2$. Understanding such a phase transition (and, more generally, threshold phenomena in randomised linear inverse problems) has been the object of formidable efforts by many researchers during the last decade, ranging from the seminal contributions by Cand\`es, Romberg and Tao \cite{crt1, crt2}, Donoho \cite{donoho1, donoho2} and Donoho and Tanner \cite{dt}, to the works of Rudelson and Vershynin \cite{RuVe08} and Ameluxen {\it et al.} \cite{AmLo13} (see \cite[Section 3]{c}, and the references therein, for a vivid description of the dense history of the subject). In particular, reference \cite{AmLo13} contains the first proof of the fundamental fact that the above described threshold phenomenon can be explained by the Gaussian concentration of the intrinsic volumes of the descent cone of the $\ell_1$ norm at ${\bf x}_0$ around its statistical dimension. In what follows, we shall further refine such a finding by showing that, for large values of $d$, the phase transition for the exact recovery of ${\bf x}_0$ has an almost exact Gaussian nature, following from the general quantitative CLTs for conic intrinsic volumes stated at Point 2 of Theorem \ref{t:mainabstract}.
\medskip

The next statement provides finite sample estimates, valid in any dimension. Note that we use the symbol $\lfloor a \rfloor$ to indicate the integer part of a real number $a$.
}
\begin{theorem}[Finite sample]\label{flag}
Let ${\bf x}_0\in\mathbb{R}^d$ and let $C$ be the descent cone of the $\ell_1$ norm $\|\cdot\|_1$ at ${\bf x}_0$. Further, let $V$ be the random variable defined by (\ref{VC}), set $\delta=E[V]$ to be the statistical dimension of $C$,  and $\tau^2 = {\rm Var}(V)$.
Let $T_{\delta,\tau}$ be the set of real numbers $t$ such that the number of observations $$m_t := \lfloor \delta+t{ \tau}\rfloor$$ lies between 1 and $d$. Fix $t\in T_{\delta,\tau}$.  Let $A_t \in \mathbb{R}^{m_t \times d}$ have independent entries, each with the standard normal ${\cal N}(0,1)$ distribution and let ${\cal F}_t=\{{\bf x}\in \mathbb{R}^d: A_t{\bf x}= A_t{\bf x}_0\}$. Consider the convex program
\[
 ({\bf CP}_t):\quad \min \|{\bf x}\|_1 \quad\mbox{subject to ${\bf x}\in {\cal F}_t$.}
\]
Then, for $\delta \ge 8$ one has the estimate
\begin{eqnarray}\label{e:iiro}
&&\sup_{t\in T_{\delta,\tau}}\left| P\left\{{\bf x}_0\mbox{ is the unique solution of $ ({\bf CP}_t)$}\right\} -
\frac{1}{\sqrt{2\pi}}\int_{-\infty}^t e^{-u^2/2}du\right| \\ &&\quad\quad\quad\quad\quad\quad\quad\quad\quad\quad\quad\quad\quad\quad\quad \leq h(\delta) + \frac{48}{\sqrt{\alpha \log^+(\alpha\sqrt{2}\delta)}}+ \frac{1}{\sqrt{2\pi\tau^2}},\notag
\end{eqnarray}
where $\alpha := \tau^2/\delta$, 
and $h(\delta)$ given by \eqref{e:xox}.
\end{theorem}

\medskip

\begin{remark}{\rm
\begin{enumerate}

\item The estimate \eqref{e:iiro} implies that, for a fixed $d$ and up to a uniform explicit error, the mapping
$$
t\mapsto P\left\{{\bf x}_0\mbox{ is the unique solution of $ ({\bf CP}_t)$}\right\},
$$
(expressing the probability of recovery as a function of $m_t$) can be approximated by the standard Gaussian distribution function $t\mapsto \Phi(t):= \frac{1}{\sqrt{2\pi}}\int_{-\infty}^t e^{-u^2/2}du$, thus demonstrating the Gaussian nature of the threshold phenomena described above. To better understand this point, fix a small $\alpha\in (0,1)$, and let $y_\alpha$ be such that $\Phi(y_\alpha) = 1-\alpha$. Then, standard computations imply that (up to the uniform error appearing in \eqref{e:iiro}) the probability $$P\left\{{\bf x}_0\mbox{ is the unique solution of $ ({\bf CP}_{y_\alpha})$}\right\}$$ is bounded from below by $1-\alpha$, whereas $P\left\{{\bf x}_0\mbox{ is the unique solution of $ ({\bf CP}_{-y_\alpha})$}\right\}$ is bounded from above by $\alpha$. Using the explicit expressions $ m_{-y_\alpha} =\lfloor \delta -y_\alpha \tau\rfloor$ and $m_{y_\alpha} = \lfloor \delta +y_\alpha \tau\rfloor$, one therefore sees that the transition from a negligible to an overwhelming probability of exact reconstruction takes place within a band of approximate length $2y_\alpha \tau \leq 2y_\alpha \sqrt{2\delta}$, centered at $\delta$. In particular, if $\delta\to \infty$, then the length of such a band becomes negligible with respect to $\delta$, thus accounting for the sharpness of the phase transition. Sufficient conditions, ensuring that $\alpha = \tau^2/\delta$ is bounded away from zero when $\delta\to\infty$, are given in Theorem \ref{t:flag2}.

\item Define the mapping $\psi : [0,1]\to [0,1]$ as
\begin{equation}\label{e:psi}
\psi(\rho) := \inf_{\gamma\geq 0} \left\{ \rho (1+\gamma^2) +(1-\rho) E[(|N|-\gamma)_+^2]\right\},
\end{equation}
where $N\sim \mathcal{N}(0,1)$. The following estimate is taken from \cite[Proposition 4.5]{AmLo13}: under the notation and assumptions of Theorem \ref{flag}, if ${\bf x_0}$ is $s$-sparse, then
\begin{equation}\label{e:uld}
\psi(s/d)-\frac{2}{\sqrt{sd}}\leq \frac{\delta}{d}\leq \psi(s/d).
\end{equation}
Moreover, as shown in \cite[Proposition 3.10]{ChRe12} one has the upper bound $\delta \leq 2s \log(d/s)+ 5s/4$, an estimate which is consistent with the classical computations contained in \cite{donoho2}.
\end{enumerate}

}
\end{remark}

}

\noindent{\it Proof of Theorem \ref{flag}}. We divide the proof into three steps.

\smallskip

\noindent{\it Step 1}. We first show that ${\bf x}_0$ is the unique solution of (${\bf CP}_t$) if and only if
$C\cap {\rm Null}(A_t)=\{{\bf 0}\}$. Indeed, assume that ${\bf x}_0$ is the unique solution of (${\bf CP}_t$) and let ${\bf y}\in C\cap {\rm Null}(A_t)$. Since ${\bf y}\in C$, there exists $\tau>0$ such that ${\bf x}:={\bf x}_0+\tau {\bf y}$ satisfies $\|{\bf x}\|_1 \leq \|{\bf x}_0\|_1$. Since ${\bf y}\in{\rm Null}(A)$ one has ${\bf x}\in {\cal F}_t$. As ${\bf x}$ is feasible the inequality $\|{\bf x}\|_1 < \|{\bf x}_0\|_1$ would contradict the assumption that ${\bf x}_0$ solves (${\bf CP}_t$). On the other hand, the equality $\|{\bf x}\|_1 = \|{\bf x}_0\|_1$ would contradict the assumption that ${\bf x}_0$ solves (${\bf CP}_t$) uniquely if ${\bf x} \not = {\bf x}_0$. Hence ${\bf y}=0$, so $C\cap {\rm Null}(A_t)=\{{\bf 0}\}$. Now assume that $C\cap {\rm Null}(A_t)=\{{\bf 0}\}$ and let ${\bf x}$ denote any solution of (${\bf CP}_t$) (note that such an ${\bf x}$ necessarily exists).
Set ${\bf y}={\bf x}-{\bf x}_0$. Of course, ${\bf y}\in{\rm Null}(A_t)$. Moreover, by definition of ${\bf x}$ and that ${\bf x}_0 \in {\cal F}_t$ one
has $\|{\bf x}\|_1=\|{\bf x}_0+{\bf y}\|_1\leq \|{\bf x}_0\|_1$, implying in turn that ${\bf y}\in C$. Hence, ${\bf y}={\bf 0}$ and ${\bf x}={\bf x}_0$, showing that ${\bf x}_0$ is the unique solution to (${\bf CP}_t$).\\

\noindent{\it Step 2}. We show\footnote{ This is a well-known result: we provide a proof for the sake of completeness.} that ${\rm Null}(A_t)\overset{\rm Law}{=}Q(\mathbb{R}^{d-m_t}\times\{{\bf 0}\})$ for $Q$ a uniformly random $d \times d$ orthogonal matrix. Both ${\rm Null}(A_t)$ and $Q(\mathbb{R}^{d-m_t}\times\{{\bf 0}\})$ belong almost surely to
the Grassmannian $G_{d-m_t}(\mathbb{R}^d)$, the set of all $(d-m_t)$-dimensional subspaces of $\mathbb{R}^d$. Defining the distance between two subspaces as the Hausdorff distance between the unit balls of those subspaces makes $G_{d-m_t}(\mathbb{R}^d)$ into a compact metric space. The metric is invariant under the action of the orthogonal group $O(d)$, and the action is transitive on $G_{d-m_t}(\mathbb{R}^d)$. Therefore, there exists a unique probability measure on $G_{d-m_t}(\mathbb{R}^d)$ that is invariant under the action of the orthogonal group. The law of the matrix $A$, having independent standard Gaussian entries, is orthogonaly invariant. Therefore,
$P({\rm Null}(A_t)\in X)=P({\rm Null}(A_t)\in R(X))$ for any $R\in O(d)$ and any measurable subset $X\subset G_{d-m_t}(\mathbb{R}^d)$.
On the other hand, it is clear that one also has $P(Q(\mathbb{R}^{d-m_t}\times\{{\bf 0}\})\in X)=P(Q(\mathbb{R}^{d-m_t}\times\{{\bf 0}\})\in R(X))$ for any $R\in O(d)$ and any measurable subset $X\subset G_{d-m_t}(\mathbb{R}^d)$.
Therefore, the claim follows by uniqueness of the probability measure on $G_{d-m_t}(\mathbb{R}^d)$ invariant under the action of $O(d)$.\\

\noindent{\it Step 3}. Combining Steps 1 and 2 we find
\[
P({\bf x}_0\mbox{ is the unique solution of (${\bf CP}_t$)})=
P(C\cap Q(\mathbb{R}^{d-m_t}  \times\{{\bf 0}\}) = \{{\bf 0}\}),
\]
where $Q$ is a uniformly random orthogonal matrix.
{ On the other hand, with $\overline{C}$ denoting the closure of $C$,
\[
P(C\cap Q(\mathbb{R}^{d-m_t}  \times\{{\bf 0}\}) = \{{\bf 0}\})=
P(\overline{C}\cap Q(\mathbb{R}^{d-m_t}  \times\{{\bf 0}\}) = \{{\bf 0}\}).
\]
As a result of this subtle point, that follows from the discussion of touching probabilities located in \cite[pp. 258--259]{sw},
we may and will assume in the rest of the proof that $C$ is closed.}
By the Crofton formula (see \cite[formula (5.10)]{AmLo13})
\bea \label{eq:Crofton}
P(C\cap Q(\mathbb{R}^{d-m_t}  \times\{{\bf 0}\}) = \{{\bf 0}\}) = 1-2h_{m_t+1}(C) \qmq{where}
h_k(C) = \sum_{j=k, j-k\,{\rm even}}^d v_j(C).
\ena
Combining \eqref{eq:Crofton} with the interlacing relation stated in \cite[Proposition 5.9]{AmLo13}, that states
\bea \label{eq:interlace}
P(V \le m_t-1) \le 1-2h_{m_t+1}(C) \le P(V \le m_t)
\ena
yields
\[
P(V \le m_t-1)\leq
P\left\{{\bf x}_0\mbox{ is the unique minimizer of (${\bf CP}_t$)}\right\}
\leq P(V \le m_t).
\]
But,
\begin{eqnarray*}
P(V \le m_t-1) &=& P\big(V \le \lfloor \delta
+t\tau\rfloor -1 \big)\\
&\geq& P\big(V\le
\delta
+t\tau - 1\big)= P\left(\frac{V-\delta}{\tau} \le t- \frac{1}{\tau} \right)
\end{eqnarray*}

and
\begin{eqnarray*}
P(V \le m_t) &=& P\big(V \le  \lfloor \delta
+t\sqrt{\tau }\rfloor \big)\\
&\leq& P\big(V \le  \delta +t\tau + 1\big)
= P\left(\frac{V-\delta}{ \tau    } \le t
+ \frac{1}{\tau}\right).
\end{eqnarray*}
 { The conclusion now follows from \eqref{e:kintro}, as well as from the fact that the standard Gaussian density on $\mathbb{R}$ is bounded by $(2\pi)^{-1/2}$.
\bbox
\medskip

The next result provides natural sufficient conditions, in order for a sequence of linear inverse problems to display exact Gaussian fluctuations in the high-dimensional limit.

\begin{theorem}[Asymptotic Gaussian phase transitions]\label{t:flag2} Let $s_n,\,d_n$, $n\geq 1$ be integer-valued sequences diverging to infinity, and assume that $s_n\leq d_n$. For every $n$, let ${\bf x}_{n,0}\in \mathbb{R}^{d_n}$
be $s_n$-sparse, denote by $C_n$ the descent cone of the $\ell_1$ norm at ${\bf x}_{n,0}$ and write $\delta_n = \delta_{C_n} = E[V_{C_n}]$ and $\tau^2_n = \tau^2_{C_n} = {\rm Var}(V_{C_n})$. For every real number $t$, write $$m_{n,t} :=
\left\{
\begin{array}{lll}
1,&&\mbox{if $\lfloor \delta_{n}+t \tau_n\rfloor<1$}\\
\lfloor \delta_{n}+t \tau_n\rfloor,&&\mbox{if $\lfloor \delta_{n}+t \tau_n\rfloor\in [1,d_n]$}\\
d_n,&&\mbox{if $\lfloor \delta_{n}+t \tau_n\rfloor>d_n$}
\end{array}
\right.
.$$ For every $n$, let $A_{n,t} \in \mathbb{R}^{m_{n,t} \times d_n}$ be a random matrix with i.i.d. ${\cal N}(0,1)$ entries, let ${\cal F}_{n,t}=\{{\bf x}\in \mathbb{R}^{d_n}: A_{n,t}{\bf x}= A_{n,t}{\bf x}_{n,0}\}$, and consider the convex program
\[
 ({\bf CP}_{n,t} ):\quad \min \|{\bf x}\|_1 \quad\mbox{subject to ${\bf x}\in {\cal F}_{n,t}$.}
\]
Assume that there exists $\rho\in (0,1)$ (independent of $n$) such that $s_n = \lfloor \rho d_n \rfloor$. Then, as $n\to\infty$, $\liminf_n \tau_n^2/\delta_n>0$, and
$$
P\left\{{\bf x}_0\mbox{ is the unique solution of $ ({\bf CP}_{n,t})$}\right\} = \frac{1}{\sqrt{2\pi}}\int_{-\infty}^t e^{-u^2/2}du + 
O\left(\frac{1}{\sqrt{\log \delta_n}}\right),
$$
where the implicit constant in the term $O\left(\frac{1}{\sqrt{\log \delta_n}}\right)$ depends uniquely on $\rho$.
\end{theorem}

\noindent{\it Proof}. In view of the estimate \eqref{e:iiro}, the conclusion will follow if we can prove the existence of a finite constant $\alpha(\rho)>0$, uniquely depending on $\rho$, such that $\tau_n^2/\delta_n\geq \alpha(\rho)$ for $n$ sufficiently large. The existence of such a $\alpha(\rho)$ is a direct consequence of the results stated in the forthcoming Proposition \ref{p:lol1}.

\bbox

\subsubsection{Second example: low-rank matrices}

Let the inner product of two $m \times n$ matrices ${\bf U}$ and ${\bf V}$ be given by
\beas
\langle {\bf U},{\bf V} \rangle = {\rm tr}({\bf U}^T {\bf V}),
\enas
and recall that, for ${\bf X} \in \mathbb{R}^{m \times n}$, the Schatten 1 (or nuclear) norm is given by
\begin{equation}\label{schatten}
\|{\bf X}\|_{S_1} = \sum_{i=1}^{\min(m,n)} \sigma_i({\bf X}),
\end{equation}
where $\sigma_1({\bf X}) \ge \cdots \ge \sigma_{\min(m,n)}({\bf X})$ are the singular values of ${\bf X}$.
Given a matrix ${\bf A} \in \mathbb{R}^{m \times np}$, partition ${\bf A}$ as $({\bf A}_1,\ldots,{\bf A}_p)$ into blocks of sizes $m \times n$, and let ${\cal A}$ be the linear map from $\mathbb{R}^{m \times n}$ to $\mathbb{R}^p$ given by
\beas
{\cal A}({\bf X}) =(\langle {\bf X},{\bf A}_1 \rangle,\cdots,\langle {\bf X},{\bf A}_p \rangle).
\enas

Now let ${\bf X}_0 \in \mathbb{R}^{m \times n}$ be a low rank matrix, and suppose that one observes
\beas
{\bf z}={\cal A}({\bf X}_0),
\enas
where the components of ${\bf A}$ are independent with distribution ${\cal N}(0,1)$.
To recover ${\bf X}_0$ we consider the convex program
\beas
{\rm min}\|{\bf X}\|_{S_1} \qmq{subject to} {\bf X} \in {\cal F},
\qmq{where}
{\cal F}=\{{\bf X}: {\cal A}({\bf X})={\bf z}\}.
\enas

As ${\cal F}$ is the affine space ${\bf X}_0 + {\rm Null}({\cal A})$, arguing as in the previous section one can show that ${\bf X}_0$ is recovered exactly if and only if $C \cap {\rm Null}({\cal A}) = \{{\bf 0}\}$ where $C={\cal D}(\|\cdot\|_{S_1},{\bf X}_0)$, the descent cone of the Schatten 1-norm at ${\bf X}_0$.

Furthermore, ${\rm Null}({\cal A})$ is a subspace of $\mathbb{R}^{m \times n}$ of dimension $nm-p$, and is rotation invariant in the sense that for any $P \subset \{(i,j): 1 \le i \le m, 1 \le j \le n\}$ of size $p$,
\beas
{\rm Null}({\cal A})= {\cal Q}(S_P)
\enas
where ${\cal Q}$ is a uniformly random orthogonal transformation on $\mathbb{R}^{m \times n}$, and
\beas
S_P = \{ {\bf X} \in \mathbb{R}^{m \times n}: X_{ij}=0 \qm{for all $(i,j) \in P$}\}.
\enas
 Now considering the natural linear mapping between $\mathbb{R}^{m \times n}$ and $\mathbb{R}^{nm}$ that preserves inner product, one may apply the Crofton formula (5.10) and proceed as for the $\ell^1$ descent cone as above in Section 1.4.3 to deduce low rank analogues of Theorems \ref{flag} and \ref{t:flag2}. In particular, for the latter we have the following result. As the Schatten 1-norm of a matrix and its transpose are equal, without loss of generality we assume that all matrices below have at least as many columns as rows.

\begin{theorem}
For every $k \in \mathbb{N}$, let $(n_k,m_k,r_k)$ be a triple of nonnegative integers depending on $k$. We assume that $n_k \rightarrow \infty$, $m_k/n_k \rightarrow \nu \in (0,1]$ and $r_k/m_k \rightarrow \rho \in (0,1)$ as $k \rightarrow \infty$, and that for every $k$ the matrix ${\bf X}(k) \in \mathbb{R}^{m_k \times n_k}$ has rank $r_k$. Let
\beas
C_k={\cal D}(\|\cdot\|_{S_1},{\bf X}(k)), \quad \delta_k=\delta(C_k) \qmq{and} \tau_k^2={\rm Var}(V_{C_k})
\enas
denote the descent cone of the Schatten 1-norm of ${\bf X}(k)$, its statistical dimension, and the the variance of its conic intrinsic volume distribution, respectively.
For every real number $t$, write
$$p_{k,t} :=
\left\{
\begin{array}{cll}
1&&\mbox{if $\lfloor \delta_{k}+t \tau_k\rfloor<1$}\\
\lfloor \delta_{k}+t \tau_k\rfloor&&\mbox{if $\lfloor \delta_{k}+t \tau_k\rfloor\in [1,m_kn_k]$}\\
m_kn_k&&\mbox{if $\lfloor \delta_{k}+t \tau_k\rfloor>m_kn_k$}
\end{array}
\right.
.$$
For every $k$, let ${\bf A}_{k,t} \in \mathbb{R}^{m_k \times n_kp_{k,t}}$ be a random matrix with i.i.d. ${\cal N}(0,1)$ entries, let ${\cal F}_{k,t}=\{{\bf X}: {\cal A}_{k,t}({\bf X})={\cal A}_{k,t}({\bf X}(k)) \}$
and consider the convex program 
\[
 ({\bf CP}_{k,t} ):\quad {\rm min}\|{\bf X}\|_{S_1} \qmq{subject to} {\bf X} \in {\cal F}_{k,t}.
\]
Then, as $k\to\infty$, $\liminf \tau_k^2/\delta_k>0$, and
$$
P\left\{{\bf X}(k)\mbox{ is the unique solution of $ ({\bf CP}_{k,t})$}\right\} = \frac{1}{\sqrt{2\pi}}\int_{-\infty}^t e^{-u^2/2}du + O\left(\frac{1}{\sqrt{\log \delta_k}}\right),
$$
where the implicit constant in the term $O\left(\frac{1}{\sqrt{\log \delta_k}}\right)$ depends uniquely on $\nu$ and $\rho$.
\end{theorem}

}

\subsection{Connections with constrained statistical inference}\label{ss:app2}

Let $C\subset \mathbb{R}^d$ be a non-trivial closed convex cone, let ${\bf g} \sim {\cal N}({\bf 0},I_d)$ and fix a vector $\bs{\mu} \in \mathbb{R}^d$. When ${\bs \mu}$ is an element of $C$ and ${\bf y} = {\bf g} + {\bs \mu}$ is regarded as a $d$-dimensional sample of observations, then the projection $\Pi_C({\bf g}+\bs{\mu})$ is  the {\it least square estimator} of ${\bs \mu}$ {\it under the convex constraint} $C$, and the norm $\|\bs{\mu}-\Pi_C({\bf g}+\bs{\mu})\|$ measures the distance between this estimator and the true value of the parameter ${\bs \mu}$; the expectation $E\|\bs{\mu}-\Pi_C({\bf g}+\bs{\mu})\|^2$ is often referred to as the $L^2$-{\it risk} of the least squares estimator.

Properties of least square estimators and associated risks have been the object of vigorous study for several decades; see e.g. \cite{BM, BVDG, Ch14, cgs, T, vdg, w, z} for a small sample. Although several results are known about the norm $\|\bs{\mu}-\Pi_C({\bf g}+\bs{\mu})\|^2$ (for instance, concerning concentration and moment estimates -- see \cite{Ch14, cgs} for recent developments), to our knowledge no normal approximation result is available for such a random variable, yet. We conjecture that our estimate \eqref{e:hairjoke} might represent a significant step in this direction. Note that, in order to make \eqref{e:hairjoke} suitable for applications, one would need explicit lower bounds on the variance of $\|\bs{\mu}-\Pi_C({\bf g}+\bs{\mu})\|^2$ for a general ${\bs \mu}$, and for the moment such estimates seem to be outside the scope of any available technique: we prefer to think of this problem as a separate issue, and leave it open for future research.

\medskip

We conclude by observing that, as explained e.g. in \cite{d,sh} and in \cite[Chapter 3]{ss}, the {\it likelihood ratio test} (LRT) for the hypotheses $H_0:{\bs \mu} = {\bf 0}$ versus $H_1 : {\bs \mu} \in C\backslash \{{\bf 0}\}$ rejects $H_0$ when the projection $\|\Pi_C({\bf y})\|^2$ of the data ${\bf y}$ on $C$ is large. In this case, our results, together with the concentration estimates from \cite{AmLo13, McTr13}, provide information on the distribution of the test statistic under the null hypothesis. Similarly, the squared projection length $\|\Pi_{C^0}({\bf y})\|^2$ onto the polar cone $C^0$
is the LRT statistic for the hypotheses $H_0:{\bs \mu} \in C$ versus $H_1 : {\bs \mu} \in \mathbb{R}^d\backslash C$.

\subsection{Plan}
The paper is organised as follows. Section 2 deals with normal approximation results for the squared distance between a Gaussian vector and a general closed convex set. Section 3 contains {total variation bounds to the normal, and}
our main CLTs for squared norms of projections onto closed convex cones, as well as for conic intrinsic volumes. { In Section 4, we derive new upper and lower bounds on the variance of conic intrinsic volumes. Section 5 is devoted to explicit Berry-Esseen bounds for intrinsic volumes distributions, whereas the Appendix in Section 6 provides a self-contained discussion of Stein's method, Poincar\'e inequalities and associated estimates on a Gaussian space.}

\section{Gaussian Projections on Closed Convex Sets: normal approximations and concentration bounds}
Let $C\subset \mathbb{R}^d$ be a closed convex set, let $\bs{\mu} \in \mathbb{R}^d$ and
let ${\bf g} \sim {\cal N}({\bf0},I_d)$ be a normal vector. In this section, we obtain a total variation bound to the normal, and a concentration inequality, for the centered squared distance between ${\bf g}+\bs{\mu}$ and $C$, that is, for
\bea \label{def:F.sigma2.gen}
F=d^2({\bf g}+\bs{\mu},C) - E[d^2({\bf g}+\bs{\mu},C)],
\ena
where $d({\bf x},C)$ is given by \eqref{def:dxC}.
We also set $\sigma^2 = {\rm Var}(d^2({\bf g}+\bs{\mu},C))={\rm Var}(F).$ It is easy to verify that $\sigma^2$ is finite for any non empty closed convex set $C$, and equals zero if and only if $C=\mathbb{R}^d$. To exclude trivialities, we call a set $C$ {\it non-trivial} if $\emptyset \subsetneq  C \subsetneq \mathbb{R}^d$.

\medskip

The following two lemmas  are the key to our main result Theorem \ref{main}: their proofs are standard, and are provided for the sake of completeness.
\begin{lemma} \label{lem:HessBound}
Let $C$ be a non empty closed convex subset of $\mathbb{R}^d$, and let $\Pi_C({\bf x})$  the metric projection onto $C$. Then, $\Pi_C$ and $I_d-\Pi_C$ are 1-Lipschitz continuous, and the Jacobian ${\rm Jac}(\Pi_C)({\bf x}) \in \mathbb{R}^{d \times d}$ exists a.e. and satisfies
\bea \label{eq:norm.deceasing}
\|(I_d-{\rm Jac}(\Pi_C)({\bf x}))^T\,{\bf y}\| \le \|{\bf y}\| \qmq{for all ${\bf y} \in \mathbb{R}^d$.}
\ena
\end{lemma}
\noindent {\em Proof:}
Since $\Pi_C$ is a projection onto a non-empty closed convex set, by \cite[p. 340]{Roc} (see also B.3 of \cite{AmLo13}), we have that
\beas
\| \Pi_C ({\bf v})-\Pi_C({\bf u})\| \leq \|{\bf v}-{\bf u}\|\qmq{for all ${\bf u},{\bf v} \in \mathbb{R}^d$,}
\enas
that is, $\Pi_C$, and hence $I_d- \Pi_C$, are
1-Lipschitz.
Bound \eqref{eq:norm.deceasing} now follows by Rademacher's theorem and the fact that, on a Hilbert space, the operator norms of a matrix and that of its transpose are the same.
\bbox

\begin{lemma} Let $C$ be a non-empty closed convex set $C \subset \mathbb{R}^d$, and let $\Pi_C({\bf x})$  be the metric projection onto $C$. Then,
\bea \label{eq:nabladsquared}
\nabla d^2({\bf x},C) = 2\left( {\bf x}- \Pi_C({\bf x}) \right), \quad {\bf x} \in \mathbb{R}^d.
\ena
\end{lemma}
\noindent {\em Proof:} Fix an arbitrary ${\bf x}_0 \in \mathbb{R}^d$, and use the shorthand notation ${\bf v}_0 :=  {\bf x_0}- \Pi_C({\bf x_0})$. Writing $\varphi({\bf u}) := d^2({\bf x}_0+{\bf u}, C) -d^2({\bf x_0}, C) -2\langle {\bf v}_0, {\bf u}\rangle$, relation \eqref{eq:nabladsquared} is equivalent to the statement that the mapping ${\bf u}\mapsto \varphi({\bf u})$ is differentiable at ${\bf u} = {\bf 0}$, and $\nabla \varphi({\bf 0}) = {\bf 0}$. To prove this statement, we show the following stronger relation: for every ${\bf u}\in \mathbb{R}^d$, one has that $ | \varphi({\bf u})| \leq \|{\bf u}\|^2$. Indeed, the inequality $\varphi({\bf u}) \leq \|{\bf u}\|^2$ follows from the fact that $d^2({\bf x}_0+{\bf u}, C) \leq \| {\bf u}+{\bf v}_0\|^2$ and
$d^2({\bf x}_0, C) = \| {\bf v}_0\|^2$. To obtain the relation $\varphi({\bf u}) \geq -\|{\bf u}\|^2$, just observe that ${\bf u} \mapsto \varphi({\bf u})$ is a convex mapping vanishing at the origin, implying that $\varphi({\bf u}) \geq - \varphi( - {\bf u}) \geq -\|\! - \!\!{\bf u}\|^2=-\|  {\bf u}\|^2 $, where the second inequality is a consequence of the estimates deduced in the first part of the proof. This yields the desired conclusion.
\bbox

\medskip

We recall that the {\it total variation distance} between the laws of two random variables $F$ and $G$ is defined as
\bea \label{def:dtv.sets}
d_{TV}(F,G)=\sup_A |P(F\in A) -P(G\in A)|,
\ena
where the supremum runs over all the Borel sets $A \subset \mathbb{R}$. It is clear from the definition that $d_{TV}(F,G)$ is invariant under affine transformations,  in the following sense: for any $a,b\in\mathbb{R}$ with $a\neq 0$, one has
$d_{TV}(aF+b,aG+b)=d_{TV}(F,G)$.
We say that $F_n$ {\it converges to $F$ in total variation} (in symbols, $F_n \stackrel{TV}{\rightarrow}F$) if $d_{TV}(F_n,F)\to 0$ as $n\to\infty$. Note that, if $F_n \ctv F$, then $F_n \claw F$, where $\claw$ denotes convergence in distribution.

\medskip

The following statement provides a total variation bound for the normal approximation of the squared distance between a Gaussian vector with arbitrary mean and a closed convex set.

\begin{theorem} \label{main}
Let $C \subset \mathbb{R}^d$ be a non trivial closed convex set, $F$ and $\sigma^2$ as in \eqref{def:F.sigma2.gen}, and $N \sim {\cal N}(0,\sigma^2)$. Then for ${\bf g} \sim {\cal N}(0,I_d)$ and $\bs{\mu} \in \mathbb{R}^d$,
\beas 
d_{{{TV}}}(F,N) \le \frac{16\sqrt{Ed^2({\bf g},C-\bs{\mu})}}{\sigma^2}.
\enas
\end{theorem}

\noindent {\em Proof:}
 As the translation of a closed convex set is closed and convex, and
\beas
d^2({\bf g}+\bs{\mu},C)=d^2({\bf g},C-\bs{\mu})
\enas
we may replace $C$ by $C-\bs{\mu}$ and assume (without loss of generality) that $\bs{\mu}={\bf 0}$. Using Lemma \ref{l:erc} and Theorem \ref{thm:dtv} in the Appendix we deduce that
\begin{equation}\label{2nd.Pbis}
d_{TV}(F,N)\leq \frac{2}{\sigma^2}\sqrt{{\rm Var}\left(\int_0^\infty e^{-t} \langle \nabla F({\bf g}), \widehat{E}(\nabla F(\widehat{\bf g}_t))\rangle dt\right)},
\end{equation}
where
\beas
{\widehat {\bf g}}_t = e^{-t} {\bf g} + \sqrt{1-e^{-2t}}{\widehat {\bf g}},
\enas
with ${\widehat {\bf g}}$ an independent copy of ${\bf g}$, and the symbols $E$ and ${\widehat E}$ denote, respectively, expectation with respect to ${\bf g}$ and ${\widehat {\bf g}}$. Set also ${\bf E}=E \otimes {\widehat E}$.
Letting ${H(\bf{g})}$ denote the integral inside the variance in \eqref{2nd.Pbis}, by \eqref{eq:nabladsquared} we have
\bea \label{eq:defG}
{H(\bf{g})}
= 4\int_0^\infty e^{-t} \langle {\bf g}- \Pi_C({\bf g}), {\widehat E} [{\widehat {\bf g}_t}-\Pi_C({\widehat {\bf g}_t})] \rangle dt.
\ena
We bound the variance of ${H(\bf{g})}$ by the Poincar\'e inequality (see Theorem \ref{thm:Poincare} in the Appendix), which states that
\begin{equation}\label{e:abb}
\mbox{Var}({H(\bf{g})}) \le  E\|\nabla {H({\bf g})}\|^2.
\end{equation}
Applying the product rule and differentiating under the integral (justified e.g. by a dominated convergence argument), using \eqref{eq:defG}, \eqref{eq:nabladsquared} and Lemma \ref{lem:HessBound} we obtain
\begin{eqnarray}
&&{\nabla H({\bf g})} = 4\int_0^\infty e^{-t} \left(I_d -{\rm Jac}(\Pi_C)({\bf g})\right)^T\, {\widehat E} [{\widehat {\bf g}}_t-\Pi_C( {\widehat {\bf g}}_t)]  dt  \label{e:maj}\\
&&\quad\quad\quad\quad\quad\quad\quad\quad\quad+ 4\int_0^\infty e^{-t} {\widehat E} [
(I_d-{\rm Jac} (\Pi_C)( {\widehat {\bf g}}_t))^T] \,\left({\bf g}- \Pi_C({\bf g}) \right) dt.\notag
\end{eqnarray}
The expectation of the squared norm of the first term on the right-hand side of \eqref{e:maj} is given by a factor of 16 multiplying
\begin{eqnarray*}
&&E\|\int_0^\infty e^{-t} \left(I_d-{\rm Jac}(\Pi_C)({\bf g})\right)^T\, {\widehat E} [ {\widehat {\bf g}}_t-\Pi_C( {\widehat {\bf g}}_t)]  dt\|^2 \\
&\le& E\int_0^\infty e^{-t} \|\left(I_d-{\rm Jac}(\Pi_C)({\bf g})\right)^T\,{\widehat E} [{\widehat {\bf g}}_t- \Pi_C( {\widehat {\bf g}}_t)]  \|^2 dt\\
&\le& E\int_0^\infty e^{-t} \|{\widehat E} [{\widehat {\bf g}}_t-\Pi_C({\widehat {\bf g}}_t)] \|^2 dt  \le {\bf E} \int_0^\infty e^{-t} \|{\widehat {\bf g}}_t-\Pi_C({\widehat {\bf g}_t})\|^2 dt \\
&=& E\int_0^\infty e^{-t} \|{\bf g}-\Pi_C({\bf g})\|^2 dt
= E\|{\bf g}-\Pi_C({\bf g})\|^2 = Ed^2({\bf g},C),
\end{eqnarray*}
where we have used the triangle inequality, Lemma \ref{lem:HessBound}, Jensen's inequality, and the fact that $\widehat{{\bf g}_t}$ has the same distribution as ${\bf g}$ for all $t$. Applying a similar chain of inequalities, it is immediate to bound the expectation of the squared norm of the second summand in \eqref{e:maj} by the same quantity. Applying \eqref{e:abb} together with the inequality $\|{\bf x}+{\bf y}\|^2 \le 2\|{\bf x}\|^2 + 2\|{\bf y}\|^2$, we therefore deduce that ${ \rm Var}(H({\bf g}))$ is bounded by $64 Ed^2({\bf g},C)$. Substituting  this bound into \eqref{2nd.Pbis} yields the desired result.
\bbox

\medskip

To conclude the section, we present a concentration bound for  random variables of the type \eqref{def:F.sigma2.gen}.
\begin{theorem} \label{thm:conc.proj.convex}
Let $C$ be a closed convex set, and $F$ given in \eqref{def:F.sigma2.gen}.
Then,
\bea \label{eq:d^2.conc.bound}
Ee^{\xi F} \le \exp \left( \frac{2 \xi^2 Ed^2({\bf g},C-\bs{\mu})}{1-2\xi} \right),
\qmq{for all $\xi < 1/2$,}
\ena
and
\begin{equation}
P(F >t) \le \exp\left( -Ed^2({\bf g},C-\bs{\mu}) h\left(\frac{t}{2 Ed^2({\bf g},C -\bs{\mu})}\right)\right) \qmq{for all $t >0$} \label{e:rez}
\end{equation}
where
\beas
h(u)=1+u-\sqrt{1+2u}.
\enas
\end{theorem}

\noindent {\em Proof:} We reduce to the case $\bs{\mu}={\bf 0}$ as in the proof of Theorem \ref{thm:main.cone}. The arguments used in the proof of Lemma 4.9 of \cite{McTr13} for convex cones work essentially  in the same way for projections on closed convex sets:  we shall therefore provide only a quick sketch of the proof, and leave the details to the reader. Similarly to \cite{McTr13}, for ${\bf g} \sim {\cal N}({\bf 0},I_d)$ we set
\beas
H({\bf g})=\xi Z \qmq{for} Z=d^2({\bf g},C)-Ed^2({\bf g},C),
\enas
and, using \eqref{eq:nabladsquared}, we deduce that
\beas
\|\nabla H({\bf g})\|^2 =4\xi^2 \|{\bf g}- \Pi_C({\bf g})\|^2 = 4 \xi^2 d^2({\bf g},C) = 4 \xi^2\left( Z+Ed^2({\bf g},C) \right).
\enas
Proceeding as in the proof of Lemma 4.9 in \cite{McTr13}, with $Ed^2({\bf g},C)$ here replacing $\delta_C$ there, yields the bound \eqref{eq:d^2.conc.bound} on the Laplace transform of $F$. Using the terminology defined in Section 2.4 of \cite{BoLuMa13}, we have therefore shown that $F$ is sub-gamma on the right tail, with variance factor $4Ed^2({\bf g},C)$ and scale parameter $2$. The conclusion now follows by the computations in that same section of \cite{BoLuMa13}.
\bbox

\bigskip

Note that the estimate \eqref{e:rez} is equivalent to the following bound: for every $t>0$
$$
P\Big(F > \sqrt{8Ed^2({\bf g}, C-{\bf \mu})t} +2t\Big) \leq e^{-t}.
$$

\begin{remark} \label{rem:polar.bear}{\rm
Let $C$ be a closed convex cone. In \cite[Lemma 4.9]{McTr13} it is proved that, for every $\xi<\frac12$,
\begin{equation}\label{e:p}
Ee^{\xi (\|\Pi_C({\bf g})\|^2 -\delta_C) } \le \exp \left( \frac{2 \xi^2 \delta_C}{1-2\xi} \right),
\end{equation}
where ${\bf g} \sim {\cal N}({\bf 0},I_d)$ and (as before) $\delta_C = E[\|\Pi_C({\bf g})\|^2]$. This estimate can be deduced by applying the general relation \eqref{eq:d^2.conc.bound} to the {polar} cone $C^0$ in the case where ${\bs \mu} = {\bf 0}$: indeed, by virtue of \eqref{eq:polar.orth.decomp} one has that
\bea \label{eq:PiC.is.d2}
\|\Pi_C({\bf x})\|^2 = d^2({\bf x},C^0),
\ena
so that \eqref{e:p} follows immediately.}
\end{remark}

\section{Steining the Steiner formula: CLTs for conic intrinsic volumes}

\subsection{Metric projections on cones}
The goal of our analysis in this subsection is to demonstrate the following variation of Theorem \ref{main}.
\begin{theorem} \label{thm:main.cone} Let $C \subset \mathbb{R}^d$ be a non-trivial closed convex cone and let
\[
F=\|\bs{\mu}-\Pi_C({\bf g}+\bs{\mu})\|^2 - m, \qmq{with} m=E[\|\bs{\mu}-\Pi_C({\bf g}+\bs{\mu})\|^2] \qmq{and} \sigma^2={\rm Var}(F).
\]
Then for every ${\bs \mu}\in \mathbb{R}^d$,
\begin{multline*}
d_{{TV}}(F,N) \le \frac{16}{\sigma^2}\left\{
\sqrt{E\| \Pi_C({\bf g}+\bs{\mu})\|^2}+2\sqrt{m}\|\bs{\mu}\|+3\|{\bs \mu}\|^2 \right\} \\
\le \frac{16}{\sigma^2}\left\{ { \sqrt{m}(1+2\|{\bs \mu}\|) + 3\|{\bs \mu}\|^2+\|{\bs \mu}\| } \right\}.
\end{multline*}
\end{theorem}
\noindent{\em Proof:} Expanding $F$ we obtain
\beas
F=\|\bs{\mu}\|^2 + \|\Pi_C({\bf g}+\bs{\mu})\|^2-2 \langle \bs{\mu}, \Pi_C({\bf g}+\bs{\mu}) \rangle - m.
\enas
The gradient of the first and last terms above are zero, while
\beas
\nabla \|\Pi_C({\bf x}+\bs{\mu})\|^2 = 2\Pi_C({\bf x}+\bs{\mu}) \qmq{and} \nabla \langle \bs{\mu}, \Pi_C({\bf x}+\bs{\mu}) \rangle = {\rm Jac}^t(\Pi_C({\bf x}+\bs{\mu})) \bs{\mu}.
\enas
the first equality following from \eqref{eq:PiC.is.d2} and \eqref{eq:nabladsquared}, the second from the definition of the Jacobian, and Lemma \ref{lem:HessBound}, showing existence. We apply \eqref{dtv2}, and hence consider
\beas
G=\int_0^\infty e^{-t} \langle \nabla F({\bf g}), \widehat{E}(\nabla F(\widehat{\bf g}_t))\rangle dt
\qmq{where}
{\widehat {\bf g}}_t = e^{-t} {\bf g} + \sqrt{1-e^{-2t}}{\widehat {\bf g}},
\enas
with ${\widehat {\bf g}}$ an independent copy of ${\bf g}$. As before, we let $E$ and ${\widehat E}$ be expectation  with respect to ${\bf g}$ and $\widehat{{\bf g}}$, respectively, and write ${\bf E}=E \otimes {\widehat E}$.

Expanding out the inner product, we obtain
\begin{multline*}
G \\
= \int_0^\infty e^{-t} \langle 2\Pi_C({\bf g}+\bs{\mu})-2{\rm Jac}^t(\Pi_C({\bf g}+\bs{\mu})) \bs{\mu}, {\widehat E} \left( 2\Pi_C({\widehat {\bf g}_t}+\bs{\mu})-2{\rm Jac}^t(\Pi_C({\widehat {\bf g}_t}+\bs{\mu})) \bs{\mu} \right)  \rangle dt\\
=4(A_1-A_2-A_3+A_4)
\end{multline*}
where
\beas
A_1 &=& \int_0^\infty e^{-t} \langle \Pi_C({\bf g}+\bs{\mu}) , {\widehat E}  \left( \Pi_C({\widehat {\bf g}_t}+\bs{\mu})\right) \rangle dt \\
A_2&=& \int_0^\infty e^{-t} \langle \Pi_C({\bf g}+\bs{\mu}) , {\widehat E}  \left(
{\rm Jac}^t(\Pi_C({\widehat {\bf g}_t}+\bs{\mu})) \bs{\mu}
\right) \rangle dt\\
A_3 &=& \int_0^\infty e^{-t} \langle
{\rm Jac}^t(\Pi_C({\bf g}+\bs{\mu})) \bs{\mu}
, {\widehat E}  \left( \Pi_C({\widehat {\bf g}_t}+\bs{\mu}) \right) \rangle dt \qm{and}\\
A_4 &=& \int_0^\infty e^{-t} \langle
{\rm Jac}^t(\Pi_C({\bf g}+\bs{\mu})) \bs{\mu}
, {\widehat E}  \left( {\rm Jac}^t(\Pi_C({\widehat {\bf g}_t}+\bs{\mu})) \bs{\mu} \right) \rangle dt.
\enas

Exploiting \eqref{dtv2},  as well as the fact that $\sigma^2 = E[G] = 4E[A_1-A_2-A_3+A_4]$, we deduce that
\bea \nonumber
d_{{ TV}}(F,N) &\le& \frac{2}{\sigma^2} E\vert \sigma^2-4\left( A_1-A_2-A_3+A_4\right)\vert \\
&\le&
\frac{8}{\sigma^2} \sum_{i=1}^4 E\vert A_i-EA_i \vert
\le \frac{8}{\sigma^2} \left( B_1+B_2+B_3+B_4\right),\label{B1234bound}
\ena
where
\beas
B_1=\sqrt{{\rm Var}(A_1)} \qmq{and} B_j=2E|A_j| \qmq{for $j=2,3,4$.}
\enas

One has that
\beas
B_j \le  2 E\left(\|\Pi_C({\bf g}+\bs{\mu})\|^2 \right)^{1/2}\|\bs{\mu}\|\leq 2(\sqrt{m}+\|\mu\|)\|\mu\| \qmq{for $j=2,3$ and} B_4 \le 2 \|\bs{\mu}\|^2,
\enas
where we have applied the Cauchy-Schwarz and triangle inequality, as well as Lemma \ref{lem:HessBound}. On the other hand, one can write
$$
A_1 =  \int_0^\infty e^{-t} \langle {\bf g}+{\bs \mu}- \Pi_{C_0}({\bf g}+\bs{\mu}) , {\widehat E}  \left( {\widehat {\bf g}}_t+{\bs \mu}- \Pi_{C_0}({\widehat {\bf g}}_t+{\bs \mu} )\right) \rangle dt,
$$
and exploit exactly the same arguments used after formula \eqref{eq:defG} (with ${\bf g}+{\bs \mu}$ and ${\widehat {\bf g}}_t+{\bs \mu}$ replacing, respectively, ${\bf g}$ and ${\widehat {\bf g}}_t$) to deduce
\beas
B_1^2 = {\rm Var}(A_1) \le 4E[\|{\bf g}+{\bs \mu} - \Pi_{C^0}({\bf g}+{\bs \mu}) \|^2] = 4E[\|\Pi_C({\bf g}+{\bs \mu}) \|^2],
\enas
thus yielding the first claim of the theorem. The second follows from observing that
{ \beas
\sqrt{E[\|\Pi_C({\bf g}+{\bs \mu}) \|^2]} \le  \sqrt{m}+\| \bs{\mu}\|,
\enas
where we have applied the triangle inequality with respect to the norm on $\mathbb{R}^d$-valued random vectors defined by the mapping $X\mapsto \sqrt{E\|X \|^2}$.          }
 \bbox

\subsection{Master Steiner formula and Main CLTs}
\label{sec:MSE&CLTs}

As anticipated in the Introduction, the aim of this section is to obtain CLTs involving the conic intrinsic volume distributions $\{\mathcal{L}(V_{C_n})\}_{n\geq 1}$ (see Section \ref{ss:intro conic}) associated with a sequence $\{C_n\}_{n\geq1}$ of closed convex cones. The strategy for achieving this goal will consist in connecting the intrinsic volume distribution of a closed convex cone $C\subset \mathbb{R}^d$ to the squared norm of the metric projection of ${\bf g} \sim {\cal N}(0,I_d)$ onto $C$.

\medskip

Our main tool will be the powerful ``Master Steiner Formula'' stated in \cite[Theorem 3.1 and Corollary 3.2]{McTr13}. Throughout the following, we use the symbol $\chi_j^2$ to indicate the chi-squared distribution with $j$ degrees of freedom, $j=0,1,2,...\,$.
	
\begin{theorem}[Master Steiner Formula, see \cite{McTr13}]\label{t:msf} Let $C\subset \mathbb{R}^d$ be a non-trivial closed convex cone, denote by $C^0$ its {polar} cone, and write $ \{v_j : j = 0,...,d\}$ to indicate the conic intrinsic volumes of $C$. Then, for every measurable mapping $f:\mathbb{R}_+^2 \rightarrow \mathbb{R}$,
\bea \label{msf}
Ef(\|\Pi_C({\bf g})\|^2,\|\Pi_{C^0}({\bf g})\|^2) = \sum_{j=0}^d E[f(Y_j,Y_{d-j}')]v_j,
\ena
where $\{Y_j,Y_j', j=0,\ldots,d\}$ stands for a collection of independent random variables such that $Y_j,Y_j' \sim \chi_j^2,j=0,1\ldots,d$.
\end{theorem}

Observe that, somewhat more compactly, we may also express \eqref{msf} as the mixture relation
\bea \label{eq:uber.master}
(\|\Pi_C({\bf g})\|^2,\|\Pi_{C^0}({\bf g})\|^2) \elaw (Y_{V_C},Y_{V_{C^0}}')
\ena
where the integer-valued random variable $V_C$  is independent of $\{Y_j,Y_j', j=0,\ldots,d\}$, and $V_{C^0} = d-V_C$. Once combined with \eqref{eq:polar.orth.decomp} and \eqref{proj.relint.j}, in the case of a polyhedral cone $C\subset \mathbb{R}^d$, relation \eqref{eq:uber.master} reinforces the intuition that, given the dimension $j$ of the face of $C$ in which lies the projection $\Pi_C({\bf g})$, the Gaussian vector ${\bf g}$ can be written as the sum of two independent Gaussian elements, with dimension $j$ and $d-j$ respectively, whose squared lengths follow the chi-squared distribution with the same respective degrees of freedom.

Fix a non-trivial closed convex cone $C\subset \mathbb{R}^d$. In order to connect the standardized limiting distributions of $\|\Pi_C({\bf g})\|^2$ and $V_C$, we use \eqref{eq:uber.master} to deduce that
\bea \label{e:xun}
\|\Pi_C({\bf g})\|^2 \elaw \sum_{i=1}^{V_C} {X_i} =  W_C + V_C, \qmq{where} W_C = \sum_{i=1}^{V_C} ({X}_i-1),
\ena
and $\{{X}_i\}_{i\geq 1}$ denotes a collection of i.i.d. $\chi_1^2$ random variables, independent of $V_C$. Since $E{X}_i=1$, we find
$E\|\Pi_C({\bf g})\|^2=E[V_C]$, and letting $G_C$ denote the squared projection length, we have
\bea \label{def:F.sigma2.gen.cone}
G_C=\|\Pi_C({\bf g})\|^2 \qmq{and} \delta_C=E[G_C].
\ena

Similarly, applying the conditional (on $ V_C$) variance formula in \eqref{e:xun} yields, with $\tau_C^2:={\rm Var}( V_C)$ and $\sigma_C^2:={\rm Var}(G_C)$, that
\bea \label{tau2.le.sigma2.-2delta}
{\rm Var} (W_C) = 2\delta_C \qmq{and} \sigma_C^2=\tau_C^2+2\delta_C,
\ena
the latter formula recovering Proposition 4.4 of  \cite{McTr13}. Standardizing both sides of the first equality in \eqref{e:xun} we therefore obtain that
\begin{equation}\label{e:x}
\frac{G_C-\delta_C}{\sigma_C} \elaw \frac{\sqrt{2\delta_C}}{\sigma_C} \frac{W_C}{\sqrt{2\delta_C}} + \frac{\tau_C}{\sigma_C} \frac{V_C-\delta_C}{\tau_C}.
\end{equation}

The following statement, that is partially a consequence of Theorem \ref{thm:main.cone}, shows that a total variation bound to the normal for the standardized projection can be expressed in terms of the mean $\delta_C$ only. We recall that $C$ is self dual when $C^0=-C$, and that in this case $\delta_C=d/2$ by \eqref{eq:VC.polar.bear}.
\begin{proposition} \label{cor:deltabound}
	We have that
	\bea \label{eq:var.delta.bounds}
	\tau_C^2 \le 2 \delta_C \qmq{and} 2 \delta_C \le \sigma_C^2 \le 4\delta_C.
	\ena
	In addition, with $G_C$ and $\delta_C$ as in \eqref{def:F.sigma2.gen.cone} and $N\sim \mathcal{N}(0,\sigma^2_C)$, one has that
	\bea \label{dtv.C}
	d_{TV}(G_C-\delta_C,N) \le \frac{16\sqrt{\delta_C}}{\sigma_C^2} \le \frac{8}{\sqrt{\delta_C}} \qmq{and, if $C$ is self dual, then} d_{TV}(F,N) \le \frac{8 \sqrt{2}}{\sqrt{d}}.
	\ena
\end{proposition}
\noindent{\em Proof:}
Theorem 4.5 of  \cite{McTr13} yields the first bound in \eqref{eq:var.delta.bounds}. The second bound in \eqref{eq:var.delta.bounds} now follows from the second relation stated in \eqref{tau2.le.sigma2.-2delta}. The first inequality in \eqref{dtv.C} follows from the first inequality of Theorem \ref{thm:main.cone} by setting $\bs{\mu}={\bf 0}$, and the remaining claims by the lower bound on $\sigma_C^2$ in  \eqref{eq:var.delta.bounds}.
\bbox

\begin{remark}{\rm The first estimate in \eqref{dtv.C} can also be directly obtained from Theorem \ref{main} by specializing it to the case $\bs{\mu}={\bf 0}$. Indeed, writing $C^0$ for the dual cone of $C$, one has that $\|\Pi_C({\bf g })\|^2 = d^2({\bf g}, C^0)$: the conclusion then follows by applying Theorem \ref{main} to the random variable $F = d^2({\bf g}, C^0) - Ed^2({\bf g}, C^0)$.}
\end{remark}

\medskip

We now consider normal limits for the conic intrinsic volumes. Explicit Berry-Esseen bounds will be presented in Theorem \ref{thm:Linfty.VC}.

\begin{theorem} \label{clt.V}
	Let $\{d_n:n\ge 1\}$ be a sequence of non-negative integers and let $\{C_n \subset \mathbb{R}^{d_n}: n \ge 1\}$ be a collection of non-trivial closed convex cones such that $\delta_{C_{n}} \rightarrow \infty$.
	For notational simplicity, write $\delta_n$, $\sigma_n$, $\tau_n$, etc., instead of $\delta_{C_{n}}$, $\sigma_{C_{n}}$, $\tau_{C_{n}}$, etc., respectively.
	Then,
	\begin{enumerate}
		\item
		\bea \label{W.to.N}
		d_{TV}\left( \frac{W_{n}}{ \sqrt{2 \delta_{n }}} , N \right)\leq \frac{2 \sigma_{n}}{\delta_{n}},\quad\mbox{{ for all $n \ge 1$}},
		\ena
		where $N\sim \mathcal{N}(0,1)$,
		and
		$$\frac{W_n}{ \sqrt{2\delta_{n} } }\ctv \mathcal{N}(0,1),\quad \mbox{as $n\to\infty$}.$$
		\item The two random variables $\frac{W_{n}}{ \sqrt{2\delta_{n}}}$ and $\frac{V_{n}-\delta_{n}}{\tau_{n}}$ are asymptotically independent in the following sense: if $\{n_k : k\geq 1\}$ is a subsequence diverging to infinity and
		\bea \label{eq:V.conv.in.dist}
		\frac{V_{n_k} -\delta_{n_k}}{\tau_{n_k}}, \quad k\geq 1,
		\ena
		converges in distribution to some random variable $Z$, then
		$$
		\left ( \frac{W_{n_k}}{ \sqrt{2\delta_{n_k}}}, \frac{V_{n_k}-\delta_{n_k}}{\tau_{n_k}}\right) \claw (N,Z),
		$$
		where $N$ has the $\mathcal{N}(0,1)$ distribution and is stochastically independent of $Z$.
		\item {If}
		\bea \label{VC.to.N}
		\frac{V_{n}-\delta_{n}}{\tau_{n}} \claw {\cal N}(0,1), \qmq{as $n \rightarrow \infty$,}
		\ena
		{then}
		\bea \label{G.to.N}
		\frac{G_{n}-\delta_{n}}{\sigma_{n}} \claw {\cal N}(0,1), \qmq{as $n \rightarrow \infty$,}
		\ena
		{and the converse implication holds if $\liminf_{n \rightarrow \infty} \tau_n^2/\delta_n>0$.}
	\end{enumerate}
\end{theorem}

\begin{remark}{\rm
		Proposition \ref{cor:deltabound} shows that, if $\delta_{n} \rightarrow \infty$,  then  \eqref{G.to.N} holds and, provided $$\liminf \tau^2_n/\delta_n>0,$$ relation \eqref{VC.to.N} also takes place by virtue of Part 3 of Theorem \ref{clt.V}. This chain of implications, which is one of the main achievements of the present paper, corresponds to the statement of Theorem \ref{t:mainabstract} in the Introduction (exception made for the Berry-Esseen bound). Results analogous to Part 3 of Theorem \ref{clt.V} (involving general mixtures of independent $\chi^2$ random variables) can be found in Dykstra \cite{dykstra}.
	}
\end{remark}

\noindent {\em Proof of Theorem \ref{clt.V}}: Throughout the proof, and when there is no risk of confusion, we drop the subscript $n$ for readability.

\noindent (\underline{Point 1}) By \cite{Lu94}, a variable $X$ with a $\Gamma(\alpha,\lambda)$ distribution satisfies
\beas
E[Xf'(X)+(\alpha-\lambda X) f(X)]=0
\enas
for all locally absolutely continuous functions $f$ for which these expectations exist.
Hence, since conditionally on $V$, $W$ has a centered chi-squared distribution with $V$ degrees of freedom, one verifies immediately that, for every Lipschitz mapping $\phi : \mathbb{R}\to \mathbb{R}$,
$$
E\left[\frac{W}{\sqrt{2\delta}}\phi\left(\frac{W}{\sqrt{2\delta}}\right)\right] = \frac{1}{\delta}E\left[(W+V) \phi'\left(\frac{W}{\sqrt{2\delta}}\right)\right].
$$
Stein's inequality (\ref{stein}) in the Appendix therefore yields that
$$
d_{TV}\left( \frac{W}{ \sqrt{2\delta}} , N \right)\leq \frac{2}{\delta} E| W +V-\delta | \leq  \frac{2}{\delta} \sqrt{2\delta + \tau^2} = \frac{2\sigma}{\delta} \le \frac{4}{\sqrt{\delta}} \rightarrow 0
$$
using \eqref{eq:var.delta.bounds}  together with the fact that $\delta \rightarrow \infty$ by assumption.

\noindent (\underline{Point 2}) Let $\eta, \, \xi$ be arbitrary real numbers. Using that the conditional distribution ${\cal L}(W|V)$ corresponds to a centered chi-squared distribution with $V$ degrees of freedom, we have
\beas
E[e^{i\eta W}|V]=\frac{e^{-i\eta V}}{(1-2 i\eta)^{V/2}} = \exp(-V(i\eta + (1/2) \log (1-2i\eta)).
\enas
Conditioning on $V$, we obtain the following expression for the joint characteristic function of ${W}/{ \sqrt{2\delta}}$ and ${(V-\delta)}/{\tau}$:
\begin{multline} \label{eq:joint.mgf}
\psi(\eta,\xi):=E\left[e^{i\eta\frac{W}{ \sqrt{2\delta}}  + i\xi \frac{V-\delta}{\tau}}\right] = E[e^{- V(i\eta/\sqrt{2\delta}+(1/2)\log(1-2i\eta/\sqrt{2\delta}))+ i\xi \frac{V-\delta}{\tau}}]\\
= e^{\delta\left [-i\eta/\sqrt{2\delta} - \frac12 \log(1-2i\eta/\sqrt{2\delta})\right]} \times  E\left[e^{  \frac{V-\delta}{\tau}\left(i\xi -i\eta\tau/\sqrt{2\delta} - \frac{\tau}{2} \log(1-2i\eta/\sqrt{2\delta})\right)} \right].
\end{multline}
As $\delta \rightarrow \infty$, one has clearly that
$$
\delta\left [-i\eta/\sqrt{2\delta} - \frac12 \log(1-2i\eta/\sqrt{2\delta})\right]\to -\eta^2/2.
$$
Moreover, since $\tau/\delta \le \sqrt{2/\delta} \rightarrow 0$ by \eqref{eq:var.delta.bounds}, we obtain as well that
$$
i\xi -i\eta\tau/\sqrt{2\delta} - \tau/2 \log(1-2i\eta/\sqrt{2\delta})\to i\xi.
$$
Hence, letting $\psi_Z$ be the characteristic function of the limiting distribution $Z$ of the sequence in \eqref{eq:V.conv.in.dist}, we infer that
$$
\psi(\eta,\xi) \rightarrow e^{-\eta^2/2}\psi_Z(\xi),
$$
thus yielding the desired conclusion.

\noindent (\underline{Point 3}) {For both implications} it is sufficient to show that, for every subsequence $n_k$, $k \ge {1}$, of {$1,2,3,\ldots$}, there exists a further subsequence $n_{k_l}, l \ge 1$, along which the claimed distributional convergence holds. By \eqref{eq:var.delta.bounds}, $0 \le \liminf \tau^2/\delta \le \limsup \tau^2/\delta  \le 2$, so for every $n_k,k \ge 0$ there exists a further subsequence $n_{k_l}, l \ge {1}$, along which $\tau^2/\delta$ converges to a limit, say $r$, in $[0,2]$. Hence, along $n_{k_l}, l \ge {1}$, we obtain
\beas \sqrt{2\delta}/{\sigma} = \sqrt{2\delta/(2\delta+ \tau^2)} \to \sqrt{\frac{2}{2+r}} \qmq{and} \tau/\sigma \to \sqrt{\frac{r}{2+r}}.
\enas

Assume first that \eqref{VC.to.N} is satisfied. Then, according to \eqref{e:x} and Point 2 in the statement, one has that
$
\frac{G-\delta}{\sigma}
$
converges in distribution along $n_{k_l},\, l \ge 1$, to $\sqrt{\frac{2}{2+r}} N + \sqrt{\frac{r}{2+r}}{Z}$, where $N$ and ${Z}$ are two independent $\mathcal{N}(0,1)$ random variables, and we conclude that \eqref{G.to.N} holds along $n_{k_l},\, l \ge {1}$.
Now assume that \eqref{G.to.N} is satisfied {and that $\liminf_{n \rightarrow \infty} \tau_n^2/\delta_n>0$; in this case, we may assume that $\tau^2/\delta$ converges to $r \in (0,2]$ along $n_{k_l}$.} {Observe that, by virtue of boundedness in $L^2$,} the family $\{\frac{V-\delta}{\tau}\}$ is tight. Consider now a further subsequence of ${n_{k_l}}$ along which $\frac{V-\delta}{\tau}$
converges in distribution to, say, $Z$. {According to Point 2 we know that the elements of the limiting pair $(N,Z)$ are independent, and by \eqref{G.to.N} the sum $\sqrt{\frac{2}{2+r}} N + \sqrt{\frac{r}{2+r}}Z$ is normal.} By Cram\'er's theorem we conclude that both $N$ and $Z$ are normally distributed, yielding the
desired conclusion.

\bbox

\bigskip

As { Table \ref{table}} below shows, Theorem \ref{t:mainabstract} yields a central limit theorem for $G_{n}$ and $V_{n}$ for the most common examples of convex cones that appear in practice. The last two rows refer to $C_A$ and $C_{BC}$, chambers of finite reflection groups acting on ${\mathbb R}^d$, which are the normal cones to the permutahedon, and signed permutahedron, respectively. For  further definitions and properties, see e.g. \cite{AmLo13, McTr13} and the references therein.


\medskip

{
	\begin{table}[h!]
	\begin{tabular}{llll}
	{\bf Cone}& {\bf Ambient} & $\delta$ & $\tau^2$ \\ \hline \hline
		Orthant & $\mathbb{R}^d$ & $\frac{1}{2}d$  & $\frac{1}{4}d$\\
		Real Positive Semi-Definite Cone & $\mathbb{R}^{n^2}$ & $\frac{1}{4}n(n+1)$ & $\simeq \left( \frac{4}{\pi^2} - \frac{1}{4}\right)n^2$\\
		$\mbox{Circ}_\alpha$ & $\mathbb{R}^d$ & $d \sin^2 \alpha$ & $\frac{1}{2}(d-2)\sin^2(2\alpha)$\\
		$C_A$ & $\mathbb{R}^d$ & $ \sum_{k=1}^d k^{-1}$ & $\sum_{k=1}^d k^{-1}(1-k^{-1})$\\
		$C_{BC}$ & $\mathbb{R}^d$ & $\sum_{k=1}^d \frac{1}{2}k^{-1}$ & $\sum_{k=1}^d \frac{1}{2}k^{-1}(1-\frac{1}{2}k^{-1})$\\
	\end{tabular}
	\caption{\it Some {common} cones}
	\label{table}
	\end{table}

\medskip
	
	\begin{remark}{\rm The first three lines of { Table \ref{table}} are { taken} from Table 6.1 of \cite{McTr13}. The means for the permutathedron and signed permutahedron are from Section D.4.\ of \cite{AmLo13}. The expressions for the variances $\tau^2$ associated with the permutathedron and signed permutahedron can be deduced as follows.
			Let
			\beas
			q(s)=\sum_{k=0}^d v_k s^k,
			\enas
			be the probability generating function of the distribution of $V=V_{C_d}$. We have
			\beas
			q'(1)=EV \qmq{and} q''(1)=EV(V-1)
			\enas
			so in particular,
			\beas
			{\rm Var}(V)=q'(1)+q''(1)-q'(1)^2 = q'(1)+\log q(s)'' \vert_{s=1}.
			\enas

			\noindent  For the permutahedron, one can use Theorem 3 of \cite[Theorem 3]{DrKl2010} (see also the first line of Table 10 of \cite{CoMo72}) to deduce that
			\beas
			q(s)=\frac{1}{d!}\prod_{k=1}^d (s+k-1) \qmq{so that} \log q(s) = -\log d! + \sum_{k=1}^d \log (s+k-1).
			\enas
			Hence,
			\beas
			EV=q'(1) = \log q(s)' \vert_{s=1}= \left( \sum_{k=1}^d \frac{1}{s+k-1}\right)_{s=1} = \sum_{k=1}^d \frac{1}{k},
			\enas
			and
			\begin{multline*}
			{\rm Var}(V)= q'(1) +\log q(s)'' \vert_{s=1} = q'(1)-\left( \sum_{k=1}^d  \frac{1}{(s+k-1)^2} \right)_{s=1} = \sum_{k=1}^d \left( \frac{1}{k} - \frac{1}{k^2} \right).
			\end{multline*}}
	\end{remark}
	The calculation for the signed permutahedron is the analogous, but now one has to use \cite[formula (3)]{BlSa98}; see also the second line of Table 10 of \cite{CoMo72}.
}

\bigskip

We conclude the section with a statement showing that the rate of convergence appearing in \eqref{W.to.N} is {often} optimal. Also, by suitably adapting the techniques introduced in \cite{exact}, one can deduce precise information about the local asymptotic behaviour of the difference $P[W_{n} /\sqrt{2\delta_{n}} \leq x] - P[N \leq x]$, where $x\in \mathbb{R}$ and $N\sim \mathcal{N}(0,1)$.

\begin{proposition}\label{p:optimality} Let the notation and assumptions of Theorem {\ref{clt.V}} prevail, and assume further that $\tau_n^2/\delta_n\to r$ for some $r\geq 0$, as $n\to\infty$. Then, for every $x\in \mathbb{R}$ one has that, as ${n}\to\infty$,
	\begin{equation}\label{e:finn}
	\frac{\delta_{n}}{\sigma_{n}}\left(P\left[\frac{W_{n}}{\sqrt{2\delta_{n}}} \leq x\right] - P[N\leq x] \right)\longrightarrow -\sqrt{\frac{2}{18+9r}} (x^2-1) \frac{e^{-x^2/2}}{\sqrt{2\pi}}.
	\end{equation}
	As a consequence, there exists a constant $c\in (0,1)$ (independent of $n$) such that, for all $n$ sufficiently large,
	\begin{equation}
	c\frac{\sigma_{n}}{\delta_{n}}\leq d_{Kol}\left(\frac{W_{n}}{\sqrt{2\delta_{n}}}, N\right)\leq d_{TV}\left(\frac{W_{n}}{\sqrt{2\delta_{n}}}, N\right).
	\end{equation}
\end{proposition}

\noindent{\it Proof.} Fix $x\in \mathbb{R}$. {It suffices to} show that, for every sequence $n_k, k\ge 1$ diverging to infinity, there exists a subsequence $n_{k_l},l\ge 1$ along which the convergence \eqref{e:finn} takes place. Let then $n_k\to\infty$ be an arbitrary divergent sequence. By $L^2$-boundedness, the collection of the laws of the random variables $
\frac{V_{n_k} -\delta_{n_k}}{\tau_{n_k}}$, $k\geq 1$ is tight, and therefore there exists a subsequence $n_{k_l}$ such that $\frac{V_{n_{k_l} }-\delta_{n_{k_l}}}{\tau_{n_{k_l}}}$ converges in distribution to some random variable $Z$. Exploiting again $L^2$-boundedness, {which additionally implies uniform integrability}, one sees immediately that $Z$ is necessarily centered. Now let $\phi_x=\phi_{h}$ denote the solution (\ref{eq:stein.solution}) {to the Stein equation \eqref{eq:stein}} for the indicator test function $h={\bf 1}_{(-\infty, x]}.$
{By (2.8) of \cite{ChGoSh10}, $\phi_x$ is Lipschitz, so as in part 1 of the proof of Theorem \eqref{clt.V}, we have
	$$
	E\left[\frac{W_n}{\sqrt{2\delta_n}}\phi\left(\frac{W_n}{\sqrt{2\delta_n}}\right)\right] = \frac{1}{\delta_n}E\left[(W_n+V_n) \phi'\left(\frac{W_n}{\sqrt{2\delta_n}}\right)\right].
	$$
	Hence, by \eqref{eq:stein}, we obtain
	\begin{eqnarray*}
		P\left[\frac{W_{n}}{\sqrt{2\delta_{n}}} \leq x\right] - P[N\leq x]  &=&
		 E\left[\phi'_x\left(\frac{W_{n}}{\sqrt{2\delta_{n}}}\right)-\frac{W_{n}}{\sqrt{2\delta_{n}}}
		\phi_x\left(\frac{W_{n}}{\sqrt{2\delta_{n}}}\right)\right]\\
		&=&
		\frac{1}{\delta_{n}} E\left[ \phi'_x\left(\frac{W_{n}}{\sqrt{2\delta_{n}}}\right) (\delta_{n} - W_{n}-V_{n})\right].
	\end{eqnarray*}
}
{Dividing both sides by $\sigma_{n}/\delta_{n}$,} one obtains
\begin{eqnarray*}
	\frac{\delta_{n}}{\sigma_{n}} \left(P\left[\frac{W_{n}}{\sqrt{2\delta_{n}}} \leq x\right] - P[N\leq x]\right)
	=E\left[ \phi'_x\left(\frac{W_{n}}{\sqrt{2\delta_{n}}}\right) \left(-\frac{\tau_{n}}{\sigma_{n}} \frac{V_{n}-\delta_{n}}{\tau_{n}} -\frac{\sqrt{2\delta_{n}}}{\sigma_{n}}\frac{W_{n}}{\sqrt{2\delta_{n}} } \right)\right].
\end{eqnarray*}

In view of {Parts 1 and 2} of Theorem \ref{clt.V}, of formula \eqref{tau2.le.sigma2.-2delta}, and of the fact that $Z$ is centered, one has, along the subsequence $n_{k_l}$, that
$$
\frac{\delta_{n}}{\sigma_{n}} \left(P\left[\frac{W_{n}}{\sqrt{2\delta_{n}}} \leq x\right] - P[N\leq x]\right)\to -\sqrt{\frac{2}{2+r}} E[ \phi'_x(N)N],
$$
where $N\sim \mathcal{N}(0,1)$. We can now use e.g. \cite[formula (2.20)]{exact} to deduce that, for every real $x$,
$$
E[ \phi'_x(N)N] = \frac{(x^2-1)}{3} \times \frac{e^{-x^2/2}}{\sqrt{2\pi}},
$$
from which the desired conclusion follows at once.

\bbox

{ In the next section, we shall prove general upper and lower bounds for the variance of conic intrinsic volumes. In particular, these results will apply to two fundamental examples that are {\it not} covered by the estimates contained in Table \ref{table}, {and that are key in convex recovery of sparse vectors and low rank matrices}: the descent cone of the $\ell_1$ norm, and of the Schatten $1$-norm.}

{

\section{Bounds on the variance of conic intrinsic volumes}\label{s:ulbounds}

\subsection{Upper and lower bounds}

Fix $d\geq 1$, let $C\subset \mathbb{R}^d$ be a closed convex cone, and let $V=V_C$ be the integer-valued random variable associated with $C$ via relation \eqref{VC}. As before, we will denote by ${\bf g}\sim \mathcal{N}(0,I_d)$ a $d$-dimensional standard Gaussian random vector. The following statement provides  useful new upper and lower bounds on the variance of $V_C$.

\begin{theorem}\label{t:lowvc} Define
\bea \label{eq:def.v.b}
v:=\| E[\Pi_C({\bf g}) ] \|^2  \qmq{and} b:=\sqrt{d\delta_C /2},
\ena
where $\delta_C$ is the statistical dimension of $C$. Then, one has the following estimates:
\begin{equation}\label{e:esti}
\frac{\min(v^2, 4b^2)}{b} \leq{\rm Var}(V_C) \leq 2v.
\end{equation}
\end{theorem}

\medskip

\begin{remark}\label{rem:remarkable}
{\rm
\begin{enumerate}
\item[(a)]In view of the orthogonal decomposition \eqref{eq:polar.orth.decomp} and of the fact that ${\bf g}$ is a centered Gaussian vector, one has that
\begin{equation}\label{e:beer}
v =  -\langle E[\Pi_C({\bf g})],E[\Pi_{C^0}({\bf g})]\rangle = \| E[\Pi_{C^0}({\bf g}) ] \|^2,
\end{equation}
where $C^0$ is the polar of $C$. Moreover, since the mapping $x\mapsto \min(x^2, 4b^2)$ is increasing on $\mathbb{R}_+$, one has also that ${\rm Var}(V_C)\geq \min(x^2, 4b^2)/b$, for every $0\leq x<v$.

\item[(b)] An elementary consequence of \eqref{e:esti} is the intuitive fact that a closed convex cone $C$ is a subspace if and only if $v=0$, that is, if and only if $\Pi_C({\bf g})$ is a centered random vector.

\end{enumerate}
}
\end{remark}

\medskip

In order to prove Theorem \ref{t:lowvc}, we need the following auxiliary result.

\begin{lemma}[Steiner form of the conic variance]\label{l:scv} {For any closed convex cone $C$,}
\beas
{\rm Var}(V_C)=-{\rm Cov}(\|\Pi_C({\bf g})\|^2,\|\Pi_{C^0}({\bf g})\|^2).
\enas
\end{lemma}
\noindent {\em Proof:} From the Master Steiner Formula \eqref{msf}, we deduce that
\begin{eqnarray*}
{\rm Cov}(\|\Pi_C({\bf g})\|^2,\|\Pi_{C^0}({\bf g})\|^2) = \sum_{j=0}^d E[Y_j Y'_{d-j}]v_j - \delta_C(d-\delta_C)
= \sum_{j=0}^d j(d-j) v_j - \delta_C(d-\delta_C),
\end{eqnarray*}
and the conclusion follows from straightforward simplifications.
\bbox

\bigskip

\noindent {\em Proof of Theorem \ref{t:lowvc}.} (Upper bound) Using \eqref{tau2.le.sigma2.-2delta}, one has that ${\rm Var}(V_C) = {\rm Var}(\|\Pi_C({\bf g})\|^2) - 2\delta_C$. Now we apply Lemma \ref{l:erc} and Theorem \ref{thm:Poincare} in the Appendix to the mapping $F({\bf g}) = \|\Pi_{C}({\bf g})\|^2 = d^2({\bf g}, C^0)$, to obtain that
$$
{\rm Var}(\|\Pi_{C}({\bf g})\|^2) \leq  \frac12 E[\| \nabla F({\bf g})\|^2] + \frac12 \| E[ \nabla F({\bf g})] \|^2 = 2\delta_C+2v,
$$
where we have used the fact that $\nabla \|\Pi_{C}({\bf g})\|^2 = 2\Pi_C({\bf g})$, {following from \eqref{eq:PiC.is.d2} and \eqref{eq:nabladsquared}}.

\smallskip

\noindent (Lower bound) For every $t>0$, define $\hat{\bf g}_t = e^{-t} {\bf g} + \sqrt{1-e^{-2t}} \hat{\bf g}$, where $\hat{\bf g}$ is an independent copy of ${\bf g}$. The crucial step is to apply relation \eqref{Poinc.Covariance} in the Appendix to the random variables $F({\bf g})=\|\Pi_C({\bf g})\|^2$ and $G({\bf g})=\|\Pi_{C^0}({\bf g})\|^2$, obtaining that, for any $a \ge 0$,
\beas
{\rm Cov}(\|\Pi_C({\bf g})\|^2,\|\Pi_{C^0}({\bf g})\|^2) \!= \!4 {\bf E} \int_0^\infty e^{-t} \langle \Pi_C({\bf g}),\Pi_{C^0}(\widehat{\bf g}_t)\rangle dt \!\le\! 4 {\bf E} \int_a^\infty e^{-t} \langle \Pi_C({\bf g}),\Pi_{C^0}(\widehat{\bf g}_t)\rangle dt,
\enas
where we have used the definition of the polar cone $C^0$ as that set that has non-positive inner product with all elements of $C$, and ${\bf E}$ indicates expectation over ${\bf g}$ and $\hat{\bf g}$. Now write
\bea \label{eq:integrate.this}
\langle \Pi_C({\bf g}),\Pi_{C^0}(\widehat{\bf g}_t)\rangle = \langle \Pi_C({\bf g}),\Pi_{C^0}(\widehat{\bf g})\rangle + \langle \Pi_C({\bf g}),\Pi_{C^0}(\widehat{\bf g}_t)-\Pi_{C^0}(\widehat{\bf g})\rangle.
\ena
For the second term, using the fact that the projection $\Pi_{C^0}({\bf x})$ is 1-Lipschitz,
\begin{multline*}
|{\bf E}\langle \Pi_C({\bf g}),\Pi_{C^0}(\widehat{\bf g}_t)-\Pi_{C^0}(\widehat{\bf g})\rangle| \le
{\bf E}\left( \|\Pi_C({\bf g})\|\, \|\Pi_{C^0}(\widehat{\bf g}_t)-\Pi_{C^0}(\widehat{\bf g})\|\right)\\
\le {\bf E}\left( \|\Pi_C({\bf g})\|\, \|\widehat{\bf g}_t-\widehat{\bf g}\|\right)
\le \sqrt{\delta(C){\bf E}\|\widehat{\bf g}_t-\widehat{\bf g}\|^2} \le \sqrt{2d\delta(C)}e^{-t}=2be^{-t},
\end{multline*}
as
\beas
{\bf E}\|\widehat{\bf g}_t-\widehat{\bf g}\|^2 ={\bf E} \|e^{-t}{\bf g} +(\sqrt{1-e^{-2t}}-1)\widehat{\bf g}\|^2 = 2\left(1- \sqrt{1-e^{-2t}} \right)d \le 2e^{-2t}d.
\enas
Now use Lemma \ref{l:scv}: multiplying \eqref{eq:integrate.this} by $e^{-t}$, integrating over $[a,\infty)$ and taking expectation yields
\beas
-{\rm Var}(V_C)\le 4{\bf E}\int_a^{\infty} e^{-t}\langle \Pi_C({\bf g}),\Pi_{C^0}(\widehat{\bf g}_t)\rangle dt \le  4e^{-a} (-v + b e^{-a}),
\enas
showing that, for every $y\in [0,1]$,
$$
{\rm Var}(V_C)\geq 4y (v - b y).
$$
The claim now follows by maximizing the mapping $y\mapsto 4y (v - b y)$ on $[0,1]$. \bbox

\medskip

In the next two sections, we shall apply the {variance bounds of Theorem \eqref{t:lowvc}} to the descent cones of the $\ell_1$ and Schatten-1 norms.

\subsection{The descent cone of the $\ell_1$ norm {at a sparse vector}} \label{subsec:ellone}
The next result provides the key for completing the proof of Theorem \ref{t:flag2}. In the body of the proofs {in this subsection and the next}, given two positive sequences $a_n,b_n$, $n\geq 1$, we shall use the notation $a_n \approx b_n$ to indicate that $a_n/b_n \to 1$, as $n\to \infty$.

\begin{proposition}\label{p:lol1} Let the assumptions and notation of Theorem \ref{t:flag2} prevail (in particular, $s_n = \lfloor \rho d_n \rfloor$ for a fixed $\rho \in (0,1)$). Then,
\begin{equation}\label{e:corinna}
\liminf_n \frac{\tau_n^2}{\delta_n}\geq \sqrt{2}  \min\left\{2\sqrt{\frac{1}{\psi(\rho)}} \,\, ;\,\,  \frac{\rho^2\gamma(\rho)^4}{\psi(\rho)^{3/2}}  \, \right\}>0,
\end{equation}
where $\psi(\rho)$ is defined in {\eqref{e:psi}} and $\gamma =\gamma(\rho)>0$ is the unique solution to the stationary equation
$$
\sqrt{\frac2\pi}\int_{\gamma}^\infty\left(\frac{u}{\gamma}-1\right)e^{-u^2/2} du  = \frac{\rho}{1-\rho}.
$$
\end{proposition}

\noindent{\it Proof}. Since the $\ell_1$ norm is invariant with respect to signed permutations, we can assume -- without loss of generality -- that {the sparse} vector ${\bf x}_{n,0}$ has the form $(x_{n,1},...,x_{n,s_n}, 0,...,0)$, $x_{n,j}>0$. Also, by virtue of the estimate \eqref{e:uld}, one has that $\delta_n \approx s_n \psi(\rho)/\rho$. Now write
$$
v_n := \|E[\Pi_{C_n}({\bf g}_n)]\|^2 = \|E[\Pi_{C^0_n}({\bf g}_n)]\|^2, \quad n\geq 1
$$
where we have used \eqref{e:beer}, and: (i) $C_n$ is the descent cone of the $\ell_1$ norm at ${\bf x}_{n,0}$, (ii) $C_n^0$ is the polar cone of $C_n$, and (iii) ${\bf g}_n = (g_1,...,g_{d_n})$ stands for a $d_n$-dimensional standard centered Gaussian vector.

 Using the lower bound in \eqref{e:esti} together with some routine simplifications, it is easily seen that relation \eqref{e:corinna} is established if one can show that
\begin{equation}\label{e:m}
\liminf_n \frac{ v_n}{s_n} \geq \gamma(\rho)^2.
\end{equation}
To accomplish this task, we first reason as in \cite[Section B.1]{AmLo13} to deduce that, for every $n$, the polar cone $C_n^0$ has the form $\bigcup_{\gamma \geq 0}\gamma \cdot \partial \|{\bf x}_{n,0}\|_1$, where {$\partial \|{\bf x}_{n,0}\|_1$ denotes the subdifferential of the $\ell^1$ norm at ${\bf x}_{n,0}$, that collection of vectors} ${\bf z} = (z_1,...,z_{d_n}) \in \mathbb{R}^{d_n}$ such that $z_1=\cdots = z_{s_n}=1$ and $|z_j|\leq 1$, for every $j=s_n+1,...,d_n$. As a consequence, for every $n$, the projection $\Pi_{C^0_n}({\bf g}_n)$ has the form
$$
\Pi_{C^0_n}({\bf g})= (\gamma_{\rho,n},...,\gamma_{\rho,n},\, \star \,,...,\,\star\,),
$$
where the symbol `$\star$' stands for entries whose exact values are immaterial for our discussion, and $\gamma_{\rho,n} > 0$ is defined as the unique random point minimising the mapping $\gamma \mapsto F_{n,\rho}(\gamma):= \sum_{i=1}^{s_n} (g_i-\gamma)^2 + \sum_{i=s_n+1}^{d_n} (|g_i|-\gamma)^2_+$ over $\mathbb{R}_+$. This shows that $v_n \geq s_n E[\gamma_{\rho,n}]^2$: as a consequence, in order to prove that \eqref{e:m} holds it suffices to check that
\begin{equation}\label{e:k}
\liminf_n E[\gamma_{\rho,n}] \geq \gamma(\rho).
\end{equation}
The key point is now that $\gamma_{\rho,n}$ is (trivially) the unique minimiser of the {\it normalised} mapping $\gamma\mapsto \frac{1}{d_n} F_{n,\rho}(\gamma)$, and also that, in view of the strong law of large numbers, for every
$\gamma\geq 0$,
\begin{equation}\label{SLLN}
\frac{1}{d_n} F_{n,\rho}(\gamma)\longrightarrow H_\rho(\gamma):= \left\{ \rho (1+\gamma^2) +(1-\rho) E[(|N|-\gamma)_+^2]\right\},\quad \mbox{as  } n\to\infty,
\end{equation}
with probability 1.

The function  $\gamma \mapsto H_\rho(\gamma)$ is minimised at the unique point $\gamma = \gamma(\rho)>0$ given in the statement, and $F_{n,\rho}(\gamma)$ is convex by (1) of Lemma C.1 of \cite{AmLo13}.
Fix $\omega\in\Omega$ and $0<\varepsilon<\gamma(\rho)$, and set
\[
D_\varepsilon=\min_{u\in\{\pm 1\}}[
H_\rho(\gamma(\rho)+\varepsilon u)-H_\rho(\gamma(\rho))].
\]
Since $\gamma(\rho)$ is the unique minimizer of $H_\rho$, one has $D_\varepsilon>0$.
From (\ref{SLLN}) we deduce the existence of
$n_0(\omega)$ large enough such that $n\geq n_0(\omega)$ implies
\[
2\max_{v\in\{0,\pm 1\}}\left|\frac{1}{d_n}F_{n,\rho}(\gamma(\rho)+\varepsilon u)-H_\rho(\gamma(\rho)+\varepsilon u)\right| <D_\varepsilon,
\]
implying in turn, by Lemma \ref{HP}, that
\[
|\gamma_{\rho,n}-\gamma(\rho)|\leq \varepsilon.
\]
That is, with probability 1,
$$
\gamma_{\rho,n}\longrightarrow \gamma(\rho)\quad \mbox{as $n\to \infty$.}
$$

Relation \eqref{e:k} now follows from a standard application of Fatou's Lemma, and the proof of \eqref{e:corinna} is therefore achieved.
\bbox
}

\subsection{The descent cone of the Schatten 1-norm {at a low rank matrix}}
{In this section we provide lower bounds on the conic variances of the descent cones of the Schatten 1-norm (see definition (\ref{schatten})) for a sequence of low rank matrices. 

For every $k \in \mathbb{N}$, let $(n,m,r)$ be a triple of nonnegative integers depending on $k$. We drop explicit dependence of $n,m$ and $r$ on $k$ for notational ease, and continue to take $m \le n$ without loss of generality. We assume that $n \rightarrow \infty,m/n \rightarrow \nu \in (0,1]$ and $r/m \rightarrow \rho \in (0,1)$ as $k \rightarrow \infty$, and that for every $k$ the matrix ${\bf X}(k) \in \mathbb{R}^{m \times n}$ has rank $r$. Let
\beas
C_k={\cal D}(\|\cdot\|_{S_1},{\bf X}(k)), \quad \delta_k=\delta(C_k) \qmq{and} \tau_k^2={\rm Var}(V_{C_k})
\enas
denote the descent cone of the Schatten 1-norm of ${\bf X}(k)$, its statistical dimension, and the the variance of its conic intrinsic volume distribution, respectively. Proposition 4.7 of \cite{AmLo13} provides that
\bea \label{eq:shatten.delta.limit}
\lim_{k \rightarrow \infty} \frac{\delta_k}{nm} = \psi(\rho,\nu),
\ena
where $\psi: [0,1]^2 \rightarrow [0,1]$ is given by
\begin{multline}
\psi(\rho,\nu) = \inf_{\gamma \ge 0} \eta(\gamma) \qmq{with} \\
\eta(\gamma) =
\left\{\rho \nu + (1-\rho \nu)
\left[ \rho(1+\gamma^2) + (1-\rho) \int_{a_-}^{a_+} (u-\gamma)_+^2 \phi_y(u) du\right] \right\}, \label{def:eta.gamma}
\end{multline}
and $y=(\nu-\rho \nu)/(1-\rho \nu)$, $a_\pm = 1 \pm \sqrt{y}$, and
\beas
\phi_y(u) = \frac{1}{\pi y u}\sqrt{(u^2-a_-^2)(a_+^2-u^2)} \quad \mbox{for $u \in [a_-,a_+]$.}
\enas
The infimum of $\eta(\gamma)$ over $[0,\infty)$ is attained at the solution $\gamma(\nu,\rho)$ to
\beas
\int_{a_- \vee \gamma}^{a_+}\left(\frac{u}{\gamma} - 1 \right) \phi_y(u)du = \frac{\rho}{1-\rho}.
\enas
It is not difficult to verify that $\gamma(\nu,\rho)>0$ for all $\nu \in (0,1],\rho \in (0,1)$.

\begin{proposition}
For the sequence of matrices ${\bf X}(k), k \in \mathbb{N}$,
\bea \label{eq:low.bd.shatten.var}
\liminf_{k \rightarrow  \infty} \frac{\tau_k^2}{\delta_k} \ge \min \left( \frac{\sqrt{2} [\rho  (1-\nu \rho) \gamma(\nu,\rho)]^2}{\psi(\rho,\nu)^{3/2}}, \frac{2^{3/2}}{\sqrt{\psi(\rho,\nu)}} \right).
\ena
\end{proposition}

\noindent{\it Proof}. By (D.8) of \cite{AmLo13}, the subdifferential of the Schatten 1-norm at ${\bf X}(k)$ is given by
\bea \label{eq:subdiff}
\partial
\|{\bf X}(k)\|_{S_1} =
\left\{
	\left[
	\begin{array}{cc}
		{\bf I}_r & 0\\
		0   & {\bf W}
	\end{array}
	\right] \in \mathbb{R}^{m \times n}: \sigma_1(W) \le 1
	\right\},
\ena
and it generates the polar $C^0$ of the descent cone, see Corollary 23.7.1 of \cite{Roc}. Closely following the proof of Proposition 4.7 of \cite{AmLo13}, and in particular the application of the Hoffman-Wielandt Theorem, see \cite{HJ90}, Corollary 7.3.8] for the second equality below, taking ${\bf G}$ to be an $m \times n$ matrix with independent ${\cal N}(0,1)$ entries, we have
\bea \nonumber
{\rm dist}({\bf G},\gamma \cdot \partial \|{\bf X}(k)\|_{S_1})^2 &=& \Bigg| \Bigg|\left[
\begin{array}{cc}
	{\bf G}_{11}-\gamma I_r & {\bf G}_{12}\\
	{\bf G}_{21}  & 0
\end{array}
\right]
\Bigg| \Bigg|_F^2 + \inf_{\sigma_1({\bf W}) \le 1}\|{\bf G}_{22}-\gamma {\bf W}\|_F^2
\\
\label{eq:Singapore2015}
&=&\Bigg| \Bigg|\left[
\begin{array}{cc}
	{\bf G}_{11}-\gamma I_r & {\bf G}_{12}\\
	{\bf G}_{21}  & 0
\end{array}
\right]
\Bigg| \Bigg|_F^2 + \sum_{i=1}^{m-r} \left(\sigma_i({\bf G}_{22})-\gamma\right)_+^2,
\ena
with $\|\cdot\|_F$ denoting the Frobenius norm and
where ${\bf G}$ is partitioned into the $2 \times 2$ block matrix $({\bf G}_{ij})_{1 \le i,j \le 2}$ formed by grouping successive rows of sizes $r$ and $m-r$, and  successive columns of sizes $r$ and $n-r$. Hence, we obtain
\bea \label{eq:proj.Shatten}
\Pi_{C_k^0}({\bf G}) = \left[
\begin{array}{cc}
\gamma_k I_r & 0 \\
0  & \gamma_k W^*
\end{array}
\right]
\ena
for some matrix $W^*$ with largest singular value at most 1, and $\gamma_k$ the minimizer of the map $\gamma \rightarrow {\rm dist}({\bf G},\gamma \cdot \partial \|{\bf X}(k)\|_{S_1})^2$ given by \eqref{eq:Singapore2015}. As the subdifferential \eqref{eq:subdiff} is a nonempty, compact, convex subset of $\mathbb{R}^{m \times n}$ that does not contain
the origin, Lemma C.1 of \cite{AmLo13} guarantees that the map is convex.

By \cite{BaSi06}, Theorem 3.6,
\beas
\frac{1}{nm}{\rm dist}^2({\bf G},\gamma\sqrt{n-r} \cdot \partial \|{\bf X}\|_{S_1}) \rightarrow_{\rm a.s.} \eta(\gamma),
\enas
where $\eta(\gamma)$ is given in \eqref{def:eta.gamma}. Reasoning as in Section \ref{subsec:ellone} (that is, using Lemma \ref{HP} followed by Fatou's lemma), we obtain
\begin{multline} \label{eq:lim.gamma.k}
\frac{\gamma_k}{\sqrt{n-r}}={\rm argmin}\left({\rm dist}^2({\bf G},\gamma\sqrt{n-r} \cdot \partial \|{\bf X}\|_{S_1})\right) \rightarrow_{\rm a.s.} \gamma(\nu,\rho) \\
\qmq{and} \liminf_{k \rightarrow \infty} \frac{E[\gamma_k]}{\sqrt{n-r}} \ge \gamma(\nu,\rho).
\end{multline}

We now invoke Theorem \ref{t:lowvc}, and make use of b) of Remark \ref{rem:remarkable}, to compute a variance lower bound in terms of
\beas
v_k=\|E[\Pi_{C_k^0}({\bf G})]\|_F^2.
\enas
The two terms in the minimum in \eqref{e:esti} give rise to the corresponding terms in \eqref{eq:low.bd.shatten.var}. By \eqref{eq:proj.Shatten},
\beas
\| \Pi_{C^0}({\bf G})\|_F \ge \sqrt{r}\gamma_k.
\enas
Squaring, taking expectation, and applying \eqref{eq:lim.gamma.k}, we find
\bea \label{eq:vk/nm>0}
\liminf_{k \rightarrow \infty} \frac{v_k}{nm} \ge \liminf_{k \rightarrow \infty} \frac{r\gamma_k^2}{nm} = \rho  (1-\nu \rho) \gamma(\nu,\rho).
\ena

Letting $b_k=\sqrt{\delta_k nm/2}$, since \eqref{eq:shatten.delta.limit} provides that $\delta_k \approx nm\psi(\rho,\nu)$, we obtain
\beas
\liminf_{k \rightarrow \infty} \frac{v_k^2}{\delta_k b_k} = \liminf_{k \rightarrow \infty} \frac{\sqrt{2} v_k^2}{(nm)^2 \psi(\rho,\nu)^{3/2}}.
\enas
Applying \eqref{eq:vk/nm>0} now yields the first term in \eqref{eq:low.bd.shatten.var}. Next, as
\beas
\liminf_{k \rightarrow \infty} \frac{4b_k}{\delta_k} =
\liminf_{k \rightarrow \infty} 2^{3/2} \sqrt{\frac{nm}{\delta_k}},
\enas
applying \eqref{eq:shatten.delta.limit} now yields the second term in \eqref{eq:low.bd.shatten.var}, completing the proof. \bbox}

\section{Bound to the normal for $V_C$}\label{ss:boundsonv}
Fix a non-trivial convex cone $C\subset\mathbb{R}^d$, and denote by $\delta_C$ and $\tau_C$, respectively, the mean and variance of its intrinsic conic distribution. The main result of the present section is Theorem \ref{thm:Linfty.VC}, providing a bound on the $L^\infty$ norm
\bea \label{eq:def.eta}
\eta = \|F-\Phi\|_\infty = \sup_{u\in \mathbb{R}}| F(u) - \Phi(u)|
\ena
of the difference between the distribution function $F(u)$ of $(V_C-\delta_C)/\tau_C$ and $\Phi(u) = P[N\leq u] $, where $N\sim \mathcal{N}(0,1)$. In the following, we set $\log^+ x = \max \left( \log x ,0\right)$.
\begin{lemma} \label{lem:Blog.plusL}
	Let $\psi_F(t)$ and $\psi_G(t)$ denote the characteristic functions of {a mean-zero distribution with variance $1$} and the standard normal distribution $\mathcal{N}(0,1)$, respectively. If
	\bea \label{eq:psi.uniform.B.a}
	\sup_{|t| \le L} |\psi_F(t)-\psi_G(t)| \le B
	\ena
	for some positive real numbers $L$ and $B$, then
	\bea \label{eq:eps.is.L1/2}
	\eta \le  B \log^+(L)  + \frac{4}{L}.
	\ena
\end{lemma}

\noindent {\em Proof:} The result holds trivially for $L<1$, so assume $L \ge 1$. Let $h_L(x)$ be the `smoothing' density function
\beas 
h_L(x)=\frac{1- \cos Lx}{\pi L x^2},
\enas
corresponding to the distribution function $H_L(x)$, { let $\Delta(x)=F(x)-G(x)$, and}  let
\beas 
\Delta_L=\Delta * H_L \qmq{and} \eta_L=\sup|\Delta_L(x)|.
\enas
By Lemma { 3.4.10} and  the proof of Lemma { 3.4.11} of \cite{dur} we have
\bea \label{eq:eta2etaL}
\eta \le 2 \eta_L + \frac{24}{\sqrt{2} \pi^{3/2} L} \qmq{and}
\eta_L \le \frac{1}{2\pi}\int_{|t| \le L} |\psi_F(t)-\psi_G(t)| \frac{dt}{|t|}.
\ena
As $\psi_F(t)$ is a characteristic function of a mean-zero distribution with variance $1$,  it is straightforward to prove that
\beas
|\psi_F(t)-1| \le \frac{t^2}{2},
\enas
so
\beas
|\psi_F(t)-\psi_G(t)| = |(\psi_F(t)-1)-(\psi_G(t)-1)|  \le t^2.
\enas
Hence for all $\epsilon \in (0,L]$
\bea \label{eq:eps.part}
\int_{|t| \le \epsilon} |\psi_F(t)-\psi_G(t)| \frac{dt}{|t|} \le \int_{|t| \le \epsilon} |t| = \epsilon^2.
\ena
By \eqref{eq:psi.uniform.B.a},
\bea \label{eq:L.part}
\int_{\epsilon < |t| \le L}|\psi_F(t)-\psi_G(t)| \frac{dt}{|t|} \le 2B \log (L/\epsilon).
\ena
Hence, by \eqref{eq:eps.part}, \eqref{eq:L.part} and \eqref{eq:eta2etaL},
\beas
\eta \le \frac{1}{\pi} \left( \epsilon^2 + 2 B\log(L/\epsilon) + \frac{24}{\sqrt{2\pi}L}\right).
\enas
As $L \ge 1$ we may choose $\epsilon=L^{-1/2}$. The conclusion now follows.
\bbox

\begin{lemma} \label{lem:L.le.tau.8}
	Let $\tau \ge 0$ and $\delta>0$ satisfy $\tau^2 \le 2\delta$. Then, the quantity
	\bea \label{eq:choose.L}
	L=\sqrt{\frac{\tau^2}{144 \delta} \log^+\left( \frac{\tau^3}{\delta}\right) }
	\qmq{satisfies}
	L\le \tau/8.
	\ena
\end{lemma}

\noindent {\em Proof:} Consider the function on $[0,\infty)$ given by
\beas
f(x) = 2 \sqrt{2} x - e^{\frac{9x^2}{4}}, \qmq{with derivative} f'(x) = 2 \sqrt{2} - \frac{9x}{2}e^{\frac{9x^2}{4}}.
\enas
Clearly, $f'(x)$ is positive at zero and decreases strictly to $-\infty$ as $x \rightarrow \infty$. Hence $f(x)$ has a global maximum value on $[0,\infty)$ achieved at the unique solution $x_0$ to the equation
\beas
xe^{\frac{9x^2}{4}}= \frac{4\sqrt{2}}{9}.
\enas
Note that
\beas
f(x_0)=2 \sqrt{2} x_0 - e^{\frac{9x_0^2}{4}} = \frac{2 \sqrt{2}}{9x_0} \left( 9x_0^2 - \frac{9}{2 \sqrt{2}} x_0 e^{\frac{9x_0^2}{4}}\right) = g(x_0) \qmq{where} g(x)=\frac{2 \sqrt{2}}{9x} \left( 9x^2 - 2\right)
\enas
and that
\beas
f'\left( \frac{\sqrt{2}}{3} \right) = \frac{\sqrt{2}}{2}\left( 4 - 3 \sqrt{e} \right) <0.
\enas
Hence $x_0 \le \sqrt{2}/3$, and since $g(x)$ is increasing in $[0,\infty)$, we have $f(x_0)=g(x_0) \le g(\sqrt{2}/3)=0$. As $f(x_0)$ is the global maximum of $f(x)$ on $[0,\infty)$ we conclude that
\bea \label{eq:another`simple'inequality!}
2 \sqrt{2}x \le e^{\frac{9x^2}{4}}.
\ena
Using $\tau^2 \le 2\delta$ and \eqref{eq:another`simple'inequality!} we obtain
\beas
\tau^3 \le 2 \sqrt{2}\delta^{3/2} \le \delta e^{\frac{9\delta}{4}} \qmq{implying} \log \left( \frac{\tau^3}{\delta} \right) \le \frac{9\delta}{4}.
\enas
The final inequality holds with $\log$ replaced by $\log^+$ since the right hand side is always non negative. The inequality so obtained provides an upper bound on $L$ in \eqref{eq:choose.L} that verifies the claim. \bbox

\medskip

In the following theorem, for notational simplicity we will write $\delta$, $\tau$ and $\sigma$ instead of $\delta_C$, $\tau_C$ and $\sigma_C$ respectively, and also set $a \vee b = \max\{a,b\}$.
\begin{theorem} \label{thm:Linfty.VC}
	The $L^\infty$ norm $\eta$ given in \eqref{eq:def.eta} satisfies
	\bea \label{eq:fu***rootlog}
	\eta \le \frac1{108} \left( \frac{\tau}{\delta^3} \vee \frac{1}{\delta^{8/3}} \right)^{\frac{3}{16}}
	\left(
	\log^+\left( \frac{\tau^3}{\delta} \right)
	\right)
	^{\frac32}
	\log^+ \left(\frac{\tau^2}{144 \delta} \log^+\left( \frac{\tau^3}{\delta}\right)  \right) + 48\sqrt{\frac{\delta}{\tau^2 \log^+\left( \frac{\tau^3}{\delta} \right)  }}.
	\ena
\end{theorem}

\begin{remark}\label{rmk:Linfty.VC}
	{\rm The estimate \eqref{e:kintro} follows immediately from \eqref{eq:fu***rootlog} and the following inequalities, valid for $\delta\geq 8$:
\begin{eqnarray*}
 \left( \frac{\tau}{\delta^3} \vee \frac{1}{\delta^{8/3}} \right)^{\frac{3}{16}} &\leq& \frac{\sqrt{2}}{\delta^{15/32}},\\ \left(\log^+\left( \frac{\tau^3}{\delta} \right)\right)^{\frac32} &\leq& (\log 2\sqrt{2\delta}) ^{3/2} \, \leq\, (\log \delta)^{3/2},\\
 \log^+ \left(\frac{\tau^2}{144 \delta} \log^+\left( \frac{\tau^3}{\delta}\right)  \right) &\leq& \log \, (\log \delta)\,\leq\, \log \delta.
\end{eqnarray*}
The above relations all follow from the bound $\tau\leq \sqrt{2\delta}$ stated in \eqref{eq:var.delta.bounds}.
}\end{remark}

\begin{remark}\label{rem:bounds.on.h}
	{\rm
		When considering a sequence of cones such that $\liminf \tau^2/\delta >0$, the right-hand side of the bound \eqref{eq:fu***rootlog} behaves like {$O \left( 1/\sqrt{\log\delta} \right)$,} thus yielding the Berry-Esseen estimate stated in Part 2 of Theorem \ref{t:mainabstract}. However, one should note that the bound \eqref{eq:fu***rootlog} covers in principle a larger spectrum of asymptotic behaviors in the {parameters $\tau^2$ and $\delta$}: in particular, in order for the right-hand side of \eqref{eq:fu***rootlog} to converge to zero, it is not necessary that the ratio $\tau^2/\delta$ is bounded away from zero}.
\end{remark}

\noindent {\em Proof:} We show Lemma \ref{lem:Blog.plusL} may be applied with $L$ as in \eqref{eq:choose.L} and
\bea \label{eq:defB}
B= 32L^3 e^{\frac{9L^2 \delta}{\tau^2}}\frac{\delta}{\tau^3}.
\ena
Let $t \in \mathbb{R}$ satisfy $|t| \le L$.
As was done in \cite{McTr13} for the Laplace transform, the implication \eqref{e:xun} of
the Steiner formula \eqref{eq:uber.master} can be applied to show that the relationship
\bea \label{eq:VGchIdent}
Ee^{it V} = Ee^{\xi_{it} G} \qmq{with} \xi_t = \frac{1}{2}\left(1 - e^{-2t} \right)
\ena
holds between the characteristic functions of $V=V_C$ and $G=\|\Pi_C({\bf g})\|^2$.
Replacing $t$ by $t/\tau$ and multiplying by $e^{-it\delta/\tau}$ in \eqref{eq:VGchIdent}
yields the following expression for the standardized characteristic function of $V$,
\bea \label{eq:VGStandIdent}
E e^{i t\left( \frac{V-\delta}{\tau}\right) }
= Ee^{\xi_{it/\tau}G}e^{-\frac{it\delta}{\tau}}.
\ena
Comparing the characteristic function of the standardized $V$ to that of the standard normal, identity \eqref{eq:VGStandIdent} and the triangle inequality yield
\begin{multline} \label{eq:three.triangle}
\vert E e^{i t\left( \frac{V-\delta}{\tau}\right) } - e^{-t^2/2} \vert = \vert Ee^{\xi_{it/\tau}G}e^{-\frac{it\delta}{\tau}} - e^{-t^2/2} \vert\\
\le \vert Ee^{\xi_{it/\tau}G} \left(
e^{ -\frac{it\delta}{\tau} }
-
e^{\left(\frac{t^2 }{\tau^2}-\xi_{it/\tau}\right)\delta}
\right) \vert
+ e^{\frac{t^2 \delta}{\tau^2}} \vert  Ee^{\xi_{it/\tau}(G-\delta)}
-  Ee^{i t/\tau (G-\delta)}  \vert \\
+  e^{\frac{t^2 \delta}{\tau^2}}  \vert Ee^{it/\tau (G-\delta)} -e^{-\sigma^2 t^2/2\tau^2}\vert.
\end{multline}
For the final term we have used \eqref{tau2.le.sigma2.-2delta}, which shows that $2\delta-\sigma^2=-\tau^2$.
For the first two terms we will make use of the inequality
\bea \label{eq:a`simple'inequality}
|e^{(a+bi)g}-e^{cig}|  \le \left(|b-c| + |a|\right)e^{|ga|}|g|,
\ena
valid for all $a,b,c,g \in \mathbb{R}$, which follows immediately by substitution from
\beas
|e^{a+bi}-e^{ci}| &=& |e^{a+bi}-e^{a+ci}+ e^{a+ci}-e^{ci}|\\
&\le& e^a |b-c| + |e^a-1| \\
&\le& e^a |b-c| + e^{|a|}-1 \\
&\le& e^{|a|}\left( |b-c| + |a|\right).
\enas

Now using \eqref{eq:VGchIdent}, implying $\vert Ee^{\xi_{it/\tau} G} \vert= \vert Ee^{(it/\tau) V} \vert \le 1$, we bound the first term in \eqref{eq:three.triangle} by
\beas
\vert  Ee^{\xi_{it/\tau}G} \vert \, \vert
e^{ -\frac{it\delta}{\tau} }
-
e^{\left(\frac{t^2 }{\tau^2}-\xi_{it/\tau}\right)\delta}  \vert
\le \vert
e^{ -\frac{it\delta}{\tau} }
-
e^{\left(\frac{t^2 }{\tau^2}-\xi_{it/\tau}\right)\delta}  \vert = \vert e^{ci}-e^{a+bi}\vert,
\enas
where we have set
\beas
a=\frac{t^2\delta}{\tau^2}-\frac{1}{2} \left(1- \cos \left( 2t/\tau \right) \right)\delta, \quad b=-\frac{1}{2} \sin(2t/\tau) \delta, \qmq{and} c=-\frac{t \delta}{\tau},
\enas
which satisfy
\beas
|a| \le \frac{2|t|^3 \delta}{3 \tau^3} \qmq{and} |b-c| \le \frac{2|t|^3\delta}{3\tau^3}.
\enas
By \eqref{eq:var.delta.bounds} of Corollary \ref{cor:deltabound} we have $\tau^2 \le 2\delta$, and in particular we may apply Lemma \ref{lem:L.le.tau.8} to yield $|t| \le L  \le \tau/8$. Now \eqref{eq:a`simple'inequality} with $g=1$ shows that the first term is bounded by
\bea \label{eq:1triangle}
\frac{4|t|^3\delta}{3 \tau^3}e^{\frac{2t^2\delta}{\tau^2}}.
\ena
Now we write the second term as
\bea \label{eq:2nd.a}
e^{\frac{t^2 \delta}{\tau^2}}E \vert
e^{\xi_{it/\tau}(G-\delta)}
- e^{it/\tau (G-\delta)}
\vert= e^{\frac{t^2 \delta}{\tau^2}} E \vert e^{(a+bi)g}-e^{cig} \vert,
\ena
where
\beas
a=\frac{1}{2} \left(1-\cos(2t/\tau) \right), \quad b=\frac{1}{2} \sin(2t/\tau), \quad c=t/\tau \qmq{and} g=G-\delta,
\enas
for which
\beas
|a| \le \min\left( \frac{|t|}{\tau},\frac{t^2}{\tau^2} \right) \qmq{and} |b-c| \le  \frac{t^2}{\tau^2}.
\enas
Applying \eqref{eq:a`simple'inequality} and the Cauchy Schwarz inequality we may bound \eqref{eq:2nd.a} as
\bea \label{eq:2nd.b}
e^{\frac{t^2 \delta}{\tau^2}}  \frac{2t^2}{\tau^2}E\left(  e^{\frac{|t|}{\tau}|G-\delta|}|G-\delta| \right) \le
\frac{2\sigma t^2}{\tau^2}e^{\frac{t^2 \delta}{\tau^2}} \sqrt{E e^{\frac{2|t|}{\tau}|G-\delta| }}.
\ena
Recalling that $\|\Pi_C({\bf x})\|^2=d^2({\bf x},C^0)$, invoking Theorem \ref{thm:conc.proj.convex} for the polar cone $C^0$ and $\bs{\mu}={\bf 0}$, for $0 \le \xi<1/2$ inequality \eqref{eq:d^2.conc.bound} yields
\beas
Ee^{\xi|G-\delta|} &=& Ee^{\xi(G-\delta)}{\bf 1}(G-\delta \ge 0) + Ee^{-\xi(G-\delta)}{\bf 1}(G-\delta <  0) \\
&\le& Ee^{\xi(G-\delta)}+ Ee^{-\xi(G-\delta)}\\
&\le& \exp \left( \frac{2 \xi^2 \delta}{1-2\xi} \right) +  \exp \left( \frac{2 \xi^2 \delta}{1+2\xi} \right) \le 2 \exp \left( \frac{2 \xi^2 \delta}{1-2\xi} \right).
\enas
Thus, applying this bound with $\xi=2|t|/\tau$, where $\xi<1/2$ by virtue of $|t| \le \tau/8$ we obtain a bound on \eqref{eq:2nd.b}, and hence on the second term of \eqref{eq:three.triangle}, of the form
\bea \label{eq:2triangle}
\frac{2\sigma t^2}{\tau^2}e^{\frac{t^2 \delta}{\tau^2}}  \sqrt{2\exp\left( \frac{8t^2 \delta}{\tau^2(1-4|t|/\tau)}\right)}
\le  2 \sqrt{2} \frac{\sigma  t^2}{\tau^2} e^{\frac{9t^2 \delta}{\tau^2}}.
\ena
For the final term, as the function $e^{it/\tau}$ has modulus 1, Theorem \ref{main} yields
\bea \label{eq:3triangle}
e^{\frac{t^2 \delta}{\tau^2}} \vert Ee^{i t/\tau(G-\delta)} - e^{- \sigma^2 t^2/2\tau^2} \vert \le 16 \frac{\sqrt{\delta_C}}{\sigma^2}e^{\frac{t^2 \delta}{\tau^2}}.
\ena
Combining the three terms \eqref{eq:1triangle}, \eqref{eq:2triangle} and \eqref{eq:3triangle}, for $|t| \le L$ we obtain
\beas
\frac{4|t|^3\delta}{3 \tau^3}e^{\frac{2t^2\delta}{\tau^2}} + 2 \sqrt{2} \frac{\sigma  t^2}{\tau^2} e^{\frac{9t^2 \delta}{\tau^2}}+ 16\frac{\sqrt{\delta_C}}{\sigma^2}e^{\frac{t^2 \delta}{\tau^2}}
\le \left( \frac{4L^3\delta}{3 \tau^3} + 2 \sqrt{2} \frac{\sigma  L^2}{\tau^2} + 16\frac{\sqrt{ \delta_C}}{\sigma^2}\right) e^{\frac{9L^2 \delta}{\tau^2}}.
\enas
From the bounds \eqref{eq:var.delta.bounds} in  Corollary \ref{cor:deltabound}, we have
\beas
\frac{4L^3\delta}{3 \tau^3} +  \frac{2 \sqrt{2} L^2 \sigma}{\tau^2} + \frac{16 \sqrt{ \delta_C}}{\sigma^2} \le
\left( \frac{4L^3}{3} + 8L^2 +16\sqrt{2} \right) \frac{\delta}{\tau^3}.
\enas
As the bound \eqref{eq:eps.is.L1/2} holds for $L<1$, we may assume $L \ge 1$, in which case $B$ as in \eqref{eq:defB} satisfies \eqref{eq:psi.uniform.B.a} when $\psi_F$ and $\psi_G$ are the characteristic functions of $(V-\delta)/\tau$ and the standard normal, respectively. Invoking Lemma \ref{lem:Blog.plusL}, the proof is completed by specializing \eqref{eq:eps.is.L1/2} to yield \eqref{eq:fu***rootlog} for the given values of $L$ and $B$.
\bbox

\section{Appendix}

\subsection{A total variation bound}

Here, we prove the total variation bound \eqref{2nd.Pbis} used in the proof of Theorem \ref{main}. We begin with a standard lemma based on Stein's method (see \cite{NoPe12}), involving the solution $\phi_h$ to the Stein equation
\bea \label{eq:stein}
\phi'_h(x)-x\phi_h(x)=h(x)-E[h(N)]
\ena
for $N \sim {\cal N}(0,1)$ and a given test functions $h$.

\begin{lemma}
If $E[F]=0$ and $E[F^2]=1$, then
	\begin{equation}\label{stein}
	d_{TV}(F,N)\leq \sup_{\phi} |E[\phi'(F)]-E[F\phi(F)]|,
	\end{equation}
	where $N\sim {\cal N}(0,1)$ and the supremum runs over all $C^1$ functions $\phi:\mathbb{R}\to\mathbb{R}$ with $\|\phi'\|_\infty\leq 2$.
\end{lemma}

\noindent {\em Proof:} For a given $h\in C^0$ taking values in $[0,1]$, by e.g.\ (2.5) of \cite{ChGoSh10}, the unique bounded solution $\phi_h(x)$ to the Stein equation \eqref{eq:stein} is given by
\bea \label{eq:stein.solution}
	\phi_h(x)=e^{x^2/2}\int_{-\infty}^x e^{-u^2/2}(h(u)-E[h(N)])du= -e^{x^2/2}\int^{\infty}_x e^{-u^2/2}(h(u)-E[h(N)])du,
\ena
where the second equality holds since
\[
\int_{\mathbb{R}} e^{-u^2/2}(h(u)-E[h(N)])du=\sqrt{2\pi} E[h(N)-E[h(N)]]=0.
\]
One can easily check that $\phi_h$ is $C^1$. Using the first equality in \eqref{eq:stein.solution} for $x<0$, and the second one for $x>0$ one obtains that $|x\phi_h(x)|\leq e^{x^2/2}\int^{\infty}_{|x|} ue^{-u^2/2}=1$. We deduce that $|\phi'_h|_\infty\leq 2$.
Recall that the total variation distance $d_{TV}(F,G)$ (as defined in \eqref{def:dtv.sets}) may also be represented as the supremum over all measurable functions $h$ taking values in $[0,1]$.
Using this fact, together with Lusin's theorem, relation \eqref{eq:stein} and the properties of the solution $\phi_h$, we infer that
\begin{multline*}
d_{TV}(F,N)=  \sup_{h:\mathbb{R}\to [0,1]} |E[h(F)]-E[h(N)]|
\\ =\sup_{h:\mathbb{R}\to [0,1],\,h\in C^0} |E[h(F)]-E[h(N)]|
\le \sup_{\phi} |E[\phi'(F)]-E[F\phi(F)]|,
\end{multline*}
as claimed. \bbox

\medskip

To make the paper as self-contained as possible, we will also prove the total variation bound \eqref{2nd.Pbis} that was applied in the proof of Theorem \ref{main}; {this result is given, { at a slightly lesser level of generality}, as Lemma 5.3 in \cite{Ch09}.}

\smallskip

{ Given $d\geq 1$, we use the symbol $\mathbb{D}^{1,2}$ to denote the Sobolev class of all mappings $f : \mathbb{R}^d \to \mathbb{R}$ that are in the closure of the set of polynomials $p : \mathbb{R}^d \to \mathbb{R}$ with respect to the norm
$$
\| p \|_{1,2} = \left (\int_{\mathbb{R}^d} p(x)^2 d\gamma(x)\right)^{1/2} + \left (\int_{\mathbb{R}^d} \|\nabla p(x)\|^2 d\gamma(x)\right)^{1/2},
$$
where $\gamma$ stands for the standard Gaussian measure on $\mathbb{R}^d$. It is not difficult to show that a sufficient condition in order for $f$ to be a member of $\mathbb{D}^{1,2}$ is that $f$ is of class $C^1$, with $f$ and its derivatives having subexponential growth at infinity. We stress that, in general, when $f$ is in $\mathbb{D}^{1,2}$ the symbol $\nabla f$ has to be interpreted in a weak sense. See e.g. \cite[Chapters 1 and 2]{NoPe12} for details on these concepts. }

\begin{theorem} \label{thm:dtv}
{ Let $H:\mathbb{R}^d\to \mathbb{R}$ be an element of $\mathbb{D}^{1,2}$}. Let ${\bf g}\sim {\cal N}(0,I_d)$ be a standard Gaussian random vector in $\mathbb{R}^d$. Let $F=H({\bf g})$ and set $m=E[F]$ and $\sigma^2={\rm Var}(F)$. Further, for $t\ge 0$, set ${\widehat {\bf g}}_t=e^{-t}{\bf g}+\sqrt{1-e^{-2t}}\widehat{\bf g}$,
where $\widehat{\bf g}$ is an independent copy of ${\bf g}$. Write $\widehat{E}$ to indicate expectation with respect to $\widehat{{\bf g}}$. Then, with $N\sim N(m,\sigma^2)$,
\begin{equation}\label{dtv}
d_{TV}(F,N)
\le \frac{2}{\sigma^2}\sqrt{{\rm Var}\left(\int_0^\infty e^{-t} \langle \nabla H({\bf g}), \widehat{E}(\nabla H(\widehat{{\bf g}}_t))\rangle dt\right)}.
\end{equation}
\end{theorem}

\noindent {\em Proof:} Without loss of generality, assume that $m=0$ and $\sigma^2=1$. The random vector
\bea \label{eq:g.in.terms.ghats}
{\bf g}_t=\sqrt{1-e^{-2t}}{\bf g} - e^{-t}\widehat{\bf g} \qmq{is an independent copy of $\widehat{\bf g}_t$, and}  {\bf g}=e^{-t}\widehat{\bf g}_t+\sqrt{1-e^{-2t}}{\bf g}_t.
\ena
{ By a standard approximation argument, it is sufficient to show the result for $H\in C^1$, with $H$ and its derivatives having subexponential growth at infinity}. Let ${\bf E}=E\otimes \widehat{E}$. If $\varphi:\mathbb{R}\to\mathbb{R}$ is $C^1$, then using the growth conditions imposed on $H$ to carry out the interchange of expectation and integration and the integration-by-parts, one has
\begin{eqnarray}
E[F\varphi(F)]&=&E[\left(H({\bf g})-H(\widehat{\bf g})\right) \varphi(H({\bf g}))] = - \int_0^\infty \frac{d}{dt} {\bf E}[H(\widehat{\bf g}_t)\varphi(H({\bf g}))]dt\notag\\
&=&
\int_0^\infty e^{-t}  {\bf E}\langle \nabla H(\widehat{\bf g}_t), {\bf g}\rangle \varphi(H({\bf g})) dt
-
\int_0^\infty \frac{e^{-2t}}{\sqrt{1-e^{-2t}}} {\bf E}\langle \nabla H(\widehat{\bf g}_t), \widehat{\bf g}\rangle \varphi(H({\bf g})) dt\notag\\
&=&\int_0^\infty \frac{e^{-t}}{\sqrt{1-e^{-2t}}} {\bf E}\langle \nabla H(\widehat{\bf g}_t), {\bf g}_t\rangle \varphi(H(e^{-t}\widehat{\bf g}_t+\sqrt{1-e^{-2t}}{\bf g}_t)) dt\notag\\
&=& \int_0^\infty e^{-t} {\bf E}\langle \nabla H(\widehat{\bf g}_t),\nabla H(e^{-t}\widehat{\bf g}_t+\sqrt{1-e^{-2t}}{\bf g}_t)\rangle \varphi'(H(e^{-t}\widehat{\bf g}_t+\sqrt{1-e^{-2t}}{\bf g}_t)) dt\notag\\
&=& E \int_0^\infty e^{-t} \langle \nabla H({\bf g}),\widehat{E}(\nabla H(\widehat{\bf g}_t))\rangle \varphi'(H({\bf g})) dt.\label{ipp2}
\end{eqnarray}
Applying identity (\ref{ipp2}) to (\ref{stein}) yields
\begin{equation}\label{dtv2}
d_{TV}(F,N)\leq 2E\left|1-\int_0^\infty e^{-t} \langle \nabla H({\bf g}), \widehat{E}(\nabla H(\widehat{\bf g}_t))\rangle dt\right|,
\end{equation}
and for $\varphi(x)=x$ yields
\bea \label{ipp}
{\rm Var}(F)=E\int_0^\infty e^{-t} \langle \nabla H({\bf g}), \widehat{E}(\nabla H(\widehat{\bf g}_t))\rangle dt.
\ena
As ${\rm Var}(F)=1$ the conclusion (\ref{dtv}), with $\sigma^2=1$, now follows by applying the Cauchy Schwarz inequality in (\ref{dtv2}).
\bbox

\bigskip

We now prove the following useful fact that was applied in the proofs of Theorem \ref{main} and Lemma \ref{l:scv}.

\begin{lemma}\label{l:erc} Let $C$ be a closed convex subset of $\mathbb{R}^d$. Then, the mapping
$$
{\bf x} \mapsto d^2({\bf x}, C)
$$
is an element of $\mathbb{D}^{1,2}$.
\end{lemma}
\noindent {\em Proof:} It is sufficient to show that $d^2(\cdot, C)$ and its derivative have sub-exponential growth at infinity. To prove this, observe that Lemma \ref{lem:HessBound} together with the triangle inequality imply that $d(\cdot , C)$ is 1-Lipschitz, so that $ d^2({\bf x}, C) \leq 2d^2({\bf 0}, C)+2\|{\bf x}\|^2$. To conclude, use \eqref{eq:nabladsquared} in order to deduce that 
$$
\|\nabla d^2({\bf x}, C)\|= 2d({\bf x}, C)\leq  2d({\bf 0}, C)+2\|{\bf x}\|.
$$
\bbox

{ A variation of the arguments leading to the proof of \eqref{ipp2} (whose details are left to the reader) yield also the following useful result.

\begin{proposition} Let $F,G\in \mathbb{D}^{1,2}$, and let the notation adopted in the statement and proof of Theorem \ref{thm:dtv} prevail. Then,
\bea \label{Poinc.Covariance}
{\rm Cov}[F({\bf g})G({\bf g})]= {\bf E} \int_0^\infty e^{-t} \langle \nabla F({\bf g}),\nabla G(\widehat{\bf g}_t)\rangle  dt.
\ena

\end{proposition}


\subsection{An improved Poincar\'e inequality}

The next result refines the classical Poincar\'e inequality, and plays a pivotal role in {Theorems \ref{main} and \ref{t:lowvc}}.

\begin{theorem}[Improved Poincar\'e inequality]\label{thm:Poincare} Fix $d\geq 1$, let $F\in \mathbb{D}^{1,2}$, and ${\bf g} = (g_1,...,g_d)\sim \mathcal{N}(0,I_d)$. Then,
$$
{\rm Var}(F({\bf g})) \leq \frac12 E[\| \nabla F({\bf g})\|^2] + \frac12 \| E[ \nabla F({\bf g})] \|^2 \leq  E[\| \nabla F({\bf g})\|^2].
$$

\end{theorem}
\noindent {\em Proof:} The quickest way to show the estimate $ {\rm Var}(F({\bf g})) \leq \frac12 E[\| \nabla F({\bf g})\|^2] + \frac12 \| E[ \nabla F({\bf g})] \|^2$ is to adopt a spectral approach. To accomplish this task, we shall use some basic  results of Gaussian analysis, whose proofs can be found e.g. in \cite[Chapter 2]{NoPe12}. Recall that, for $k=0,1,2,...$, the {$k^{th}$} {\it Wiener chaos} associated with ${\bf g}$, written $C_k$, is the subspace spanned by all random variables of the form $\prod_{i=1}^m H_{k_i}(g_{j_i})$, where $\{H_k : k=0,1,...\}$ denotes the collection of Hermite polynomials on the real line, $k_1+\cdots + k_m = k$, and the indices $k_1,\cdots, k_m$ are pairwise distinct. It is easily checked that Wiener chaoses of different orders are orthogonal in $L^2(\Omega)$, and also that every square-integrable random variable of the type $F({\bf g})$ can be decomposed as an infinite sum of the type $F({\bf g})= \sum_{k=0}^\infty F_k({\bf g})$, where the series converges in $L^2(\Omega)$ and where, for every $k$, $F_k({\bf g})$ denotes the projection of $F({\bf g})$ on $C_k$ (in particular, $F_0({\bf g}) = E[F({\bf g})]$). This {decomposition yields in particular} that
$$
{\rm Var}(F({\bf g})) = \sum_{k=1}^\infty E[F^2_k({\bf g})].
$$
The key point is now that, if $F\in \mathbb{D}^{1,2}$, then one has the additional relations
$$
E[\| \nabla F({\bf g})\|^2] = \sum_{k=1}^\infty k E[F^2_k({\bf g})]
$$
(see e.g. \cite[Exercice 2.7.9]{NoPe12}) and
$$
E[F_1^2({\bf g})] = \| E[ \nabla F({\bf g})] \|^2,
$$
the last identity being justified as follows: if $F$ is a smooth mapping, then the projection of $F({\bf g})$ on $C_1$ is given by
$$
F_1({\bf g}) = \sum^d_{i=1} E[F({\bf g})g_i] g_i = \sum^d_{i=1}E\left[\frac{\partial F}{\partial x_i} ({\bf g})\right]g_i,
$$
and the result for a general $F\in \mathbb{D}^{1,2}$ is deduced by an approximation argument. The previous relations imply therefore that
$$
{\rm Var}(F({\bf g})) \!=\! \sum_{k=1}^\infty E[F^2_k({\bf g})]\leq E[F_1^2({\bf g})] \!+\!  \sum_{k=2}^\infty \frac k2 E[F^2_k({\bf g})] \!=\! \! \frac12 \| E[ \nabla F({\bf g})] \|^2+\frac12 E[\| \nabla F({\bf g})\|^2] \!\, .
$$
The proof is concluded by observing that, in view of Jensen inequality, $\| E[ \nabla F({\bf g})] \|^2 \leq  E[\| \nabla F({\bf g})\|^2].$
\bbox

}

\subsection{A bound on the distance to the minimizer of a convex function.}

Following an idea introduced by Hjort and Pollard \cite{HP}, one has the following lemma, providing
a bound on the distance to the minimizer of a convex function in terms of another, not necessarily convex, function.

\begin{lemma}\label{HP}
Suppose $f:[0,\infty)\to\mathbb{R}$ is a convex function, and let $g:[0,\infty)\to\mathbb{R}$ be any function. If $x_0$ is a minimizer of $f$, $y_0\in(0,\infty)$ and $\varepsilon\in(0,y_0)$, then
\begin{equation}\label{7.36}
2\max_{v\in\{0,\pm 1\}}|g(y_0+\varepsilon v)-f(y_0+\varepsilon v)| <\min_{u\in\{\pm 1\}} [g(y_0+\varepsilon u)-g(y_0)]
\end{equation}
implies $|x_0-y_0|\leq \varepsilon$.
\end{lemma}
{\it Proof}. Suppose $a:=|x_0-y_0|> \varepsilon>0$. Set $u=a^{-1}(x_0-y_0)$. Then $u\in\{\pm 1\}$, $x_0=y_0+au$ and the convexity of $f$ implies
\[
(1-\varepsilon/a)f(y_0)+(\varepsilon/a)f(x_0)\geq f(y_0+\varepsilon u).
\]
Hence
\begin{eqnarray*}
&&\frac{\varepsilon}{a}(f(x_0)-f(y_0))
\geq f(y_0+\varepsilon u)-f(y_0)\\
&=&g(y_0+\varepsilon u) - g(y_0) + [f(y_0+\varepsilon u)-g(y_0+\varepsilon u)] + [g(y_0)-f(y_0)]\\
&\geq&\min_{u\in\{\pm 1\}}  [ g(y_0+\varepsilon u)-g(y_0)] - 2\max_{v\in\{0,\pm 1\}}|g(y_0+\varepsilon v)-f(y_0+\varepsilon v)|.
\end{eqnarray*}
If (\ref{7.36}) is satisfied, then $\frac{\varepsilon}{a}(f(x_0)-f(y_0))>0$. But this contradicts that $x_0$ is a minimizer of $f$. Hence, $|x_0-y_0|> \varepsilon$ is impossible.\bbox

\bigskip
\bigskip

\noindent {\bf Acknowledgments:} The authors would like to thank John Pike for his assistance with the final two rows of Table \ref{table}, and the associated computations.

\end{document}